\documentclass[a4paper,fleqn]{cas-sc}

\usepackage[utf8]{inputenc}
\usepackage{natbib}
\usepackage{hyperref}
\usepackage{placeins}
\definecolor{brickred}{rgb}{0.8, 0.0, 0.0}
\definecolor{cyan_custom}{rgb}{0.0, 0.50, 0.67}
\usepackage{amsmath}
\usepackage{amsthm}
\usepackage{algorithm}
\usepackage{algorithmic}
\usepackage{lineno} 
\newtheoremstyle{mydef}
  {10pt} 
  {10pt} 
  {\normalfont} 
  {} 
  {\bfseries} 
  {:} 
  {.5em} 
  {} 
\theoremstyle{mydef}
\newtheorem{definition}{Definition}
\def\tsc#1{\csdef{#1}{\textsc{\lowercase{#1}}\xspace}}
\tsc{WGM}
\tsc{QE}
\tsc{EP}
\tsc{PMS}
\tsc{BEC}
\tsc{DE}

\begin{document}
\hypersetup{allcolors = cyan_custom}
\let\WriteBookmarks\relax
\def\floatpagepagefraction{1}
\def\textpagefraction{.001}



\title [mode = title]{Assessing How Ride-hailing Rebalancing Strategies Improve the Resilience of Multi-modal Transportation Systems}

\author{Euntak Lee}[]
\cormark[1] 
\ead{euntak.lee@univ-eiffel.fr}  
\credit{Conceptualization of this study, Methodology, Software}  
\affiliation{organization={Univ. Eiffel}, addressline={ENTPE, LICIT-ECO7}, postcode={F-69518}, city={Lyon}, country={France}}  
\author{Rim Slama}
\ead{rim.slamasalmi@entpe.fr}
\author{Ludovic Leclercq}
\ead{ludovic.leclercq@univ-eiffel.fr}
\credit{Data curation, Writing - Original draft preparation}
\cortext[cor1]{ Corresponding author}

\begin{abstract}
The global ride-hailing (RH) industry plays an essential role in multi-modal transportation systems by improving user mobility, particularly as first- and last-mile solutions. However, the flexibility of on-demand mobility services can lead to local supply-demand imbalances. While many RH rebalancing studies focus on nominal scenarios with regular demand patterns, it is crucial to consider disruptions—such as train line interruptions—that negatively impact operational efficiency, resulting in longer travel times, higher costs, increased transfers, and service delays. This study examines how RH rebalancing strategies can strengthen the resilience of multi-modal transportation systems against such disruptions. We incorporate RH services into systems where users choose and switch transportation modes based on their preferences, accounting for uncertainties in demand predictions that reflect discrepancies between forecasts and actual conditions. To address the stochastic supply-demand dynamics in large-scale networks, we propose a multi-agent reinforcement learning (MARL) strategy, specifically utilizing a multi-agent deep deterministic policy gradient (MADDPG) approach. The proposed framework is particularly well-suited for this problem due to its ability to handle continuous action spaces, which are prevalent in real-world transportation systems, and its capacity to enable effective coordination among multiple agents operating in dynamic and decentralized environments. Through a 900 $\text{km}^2$ multi-modal traffic simulation, we evaluate the proposed model's performance against four existing RH rebalancing strategies, focusing on its ability to enhance system resilience. The results demonstrate significant improvements in key performance indicators, including user waiting time, resilience metrics, total travel time, and travel distance.
\end{abstract}

\begin{keywords}
Ride-hailing services
\sep Rebalancing strategy 
\sep Resilience
\sep Multi-modal transportation systems
\sep Reinforcement learning
\end{keywords}

\maketitle
\section{Introduction}

Ride-hailing (RH) mobility services are expected to transform the future of multi-modal transportation systems. The global RH market is expected to grow from USD 243 billion in 2023 to USD 519 billion by 2030, at an annual growth rate of 11.44\% during the forecast period \citep{RH_Market}. RH services account for about 15\% of all trips in large cities, and users are expected to have the potential to reduce their dependence on private cars \citep{erhardt2019transportation}. Furthermore, a greater use of RH services is associated with a greater frequency of public transportation systems, such as metros, buses, and trains (\citealp{feigon2016shared}; \citeyear{feigon2018broadening}). 

The role of RH services in first and last mile trips is essential to encourage the shift from private cars to public transportation-based multi-modal trips \citep{sunitiyoso2022role}. As the RH industry continues to expand, its role in enhancing the efficiency of multi-modal transportation systems becomes increasingly crucial. Effective integration can significantly improve user mobility. Yet, alongside these opportunities lies a pressing challenge: the disruptions that frequently compromise the performance and reliability of such interconnected networks.

Disruptions in multi-modal transportation networks, whether stemming from supply or demand-side issues, pose a significant threat to system performance and user satisfaction. Supply disruptions lead to reduced availability of the road network, which can occur due to road accidents or lack of mobility services. On the other hand, demand disruptions are related to increases in transportation demand, such as during large sports events, strikes, or the arrival of mass public transportation. Excess demand results in interruptions in system functioning due to limited supply capacity. For example, public transportation routes can temporarily close, resulting in the discharge of numerous passengers at some stations. As a cascading effect, a surge in demand for limited alternative modes could cause all users to experience longer travel times, higher costs, increased transfers, and more service delays (\citealp{cadarso2013recovery}; \citealp{shahparvari2016enhancing}).

For operational efficiency, a more resilient system must be established that can withstand and absorb the negative impacts of disruptions, maintain a certain level of service, and recover quickly \citep{ju2022multilayer}. Transportation modes in multi-modal transportation systems are closely interconnected in terms of passenger flow and facilities. If disruptions affect one transportation mode, stranded passengers can opt for alternative transportation options, sequentially increasing the burden of other transportation modes. Therefore, it is vital to comprehensively consider the interrelationships among these transportation modes and develop resilient transportation systems. As part of the multi-modal transportation system, appropriate RH service operations can provide a promising way to improve resilience.

\subsection{Resilience of Transportation Systems}
Research on the resilience of transportation systems can be categorized into three areas: (i) resilience response, (ii) resilience enhancement, and (iii) resilience assessment \citep{zhang2024analysis}. First, the resilience response proposes effective reactions to disruptions. Throughout the periods, various approaches for resilience response have been suggested, such as emergency rescue capacity for the urban subway network \citep{chen2022resilience}, bus bridging services \citep{chen2022resilience2}, vehicle routing \citep{wohlgemuth2012dynamic}, traffic signal control \citep{aslani2017adaptive}, and train traffic control \citep{shakibayifar2019intelligent}. Second, the resilience enhancement focuses on strengthening a system's ability itself to withstand solidly and recover swiftly from disruptions. The proposals for securing network spare capacity include managing alternative travel diversity \citep{xu2018transportation}, providing on-demand transport services \citep{xu2022network}, and combining transportation with power networks for electric vehicles \citep{hussain2022resilience}. Lastly, the resilience assessment focuses on proposing metrics and evaluation frameworks to define the resilience of systems.

Resilience assessment has been focused on the practical performance of systems in different phases of disruption. For example, \citet{jenelius2015value} measured the robustness and redundancy of a public transport system when disruptions occur such as service rate reductions and link breakdowns. They proposed the measurements by applying the passenger welfare, or total utility, of a generalized cost function that includes in-vehicle time, waiting time, walking time, and number of transfers. \citet{xu2021optimizing} proposed a tolerance function that expects worst-case scenarios to evaluate a system's performance under rail transit disruptions. By applying the function, they suggested the optimal disruption mitigation strategy, such as platform downtime protection and bus-bridging service planning. Similarly, various metrics for different concepts have been proposed, including recovery cost \citep{henry2012generic}, centrality \citep{sun2016measuring}, minimum system performance \citep{goldbeck2019resilience}, readiness \citep{ahmadian2020quantitative}, route diversity \citep{jansuwan2021analysis}, and rapidity \citep{fang2022resilience}. However, most of the existing metrics are associated with any of the following limitations: (1) resilience from either but not both perspectives of network structure or passenger demand is considered; (2) static attributes of systems are mainly considered, neglecting dynamic features such as real-time flows of passengers and vehicles; or (3) impractical assumptions are applied as responses to disruptions are not considered. To overcome this, \citet{zhang2024analysis} recently proposed a combination of resilience metrics in four aspects: vulnerability, adaptability, robustness, and recoverability. This metric considers the impacts of disruptions on passenger choices of paths and a dynamic simulation procedure in the overall disruption phases. As it only reflects supply disruptions, metrics need to be developed that consider demand disruptions due to special situations such as COVID-19, big sport events, concerts, and departure and arrival of ships, boats, and airplanes. These aforementioned studies have been conducted within public transportation systems, establishing a foundation for developing resilient transportation systems. Multifaceted evaluation can provide insights into how the system should be assessed, which reactions should be needed to handle disruptions, and how to improve resilience. 

Efforts have been made to integrate RH mobility services into multi-modal transportation systems, demonstrating their crucial role as a supplementary mode (\citealp{alonso2018potential}; \citealp{liu2023integrating}; \citealp{cortina2023fostering}). These efforts primarily focused on improving the operational efficiency of RH services under regular scenarios with consistent demand patterns. However, to the best of our knowledge, the resilience of RH services within multi-modal transportation systems remains unexamined. For reliable operations, it is essential to explore how RH services contribute to the resilience of multi-modal transportation systems against disruptions. This study seeks to bridge this gap by studying resilience assessment while integrating RH services within multi-modal transportation systems, providing critical insights into their role in enhancing system robustness and adaptability to disruptions.

\subsection{Ride-hailing Rebalancing Strategies}
RH services offer the convenience of allowing users to request rides to any location of their choice. Leveraging these strengths, RH companies such as Uber, Lyft, and Bolt have transformed urban mobility by providing flexible, on-demand transportation options. The RH industries continues to expand rapidly; for instance, Uber reported the number of active vehicles to approximately 28,000 in Paris, France, in 2019 \citep{paris_court_2019}. However, this freedom may cause supply-demand imbalances in local areas. Users cannot request rides if there are no RH vehicles nearby, while other areas are oversupplied with vehicles. Consequently, vehicle fleet management has gained interest for efficient operations. However, operating an appropriate strategy is not a simple issue since inefficient operations can instead lead to longer vacant times and longer waiting times for drivers and passengers, respectively \citep{beojone2021inefficiency}. As extra RH vehicles cannot be deployed whenever there is a shortage, managing this imbalance is one of the critical issues \citep{lee2020optimal}.

\begin{figure}[h]
    \centering
    \includegraphics[width=0.8\textwidth]{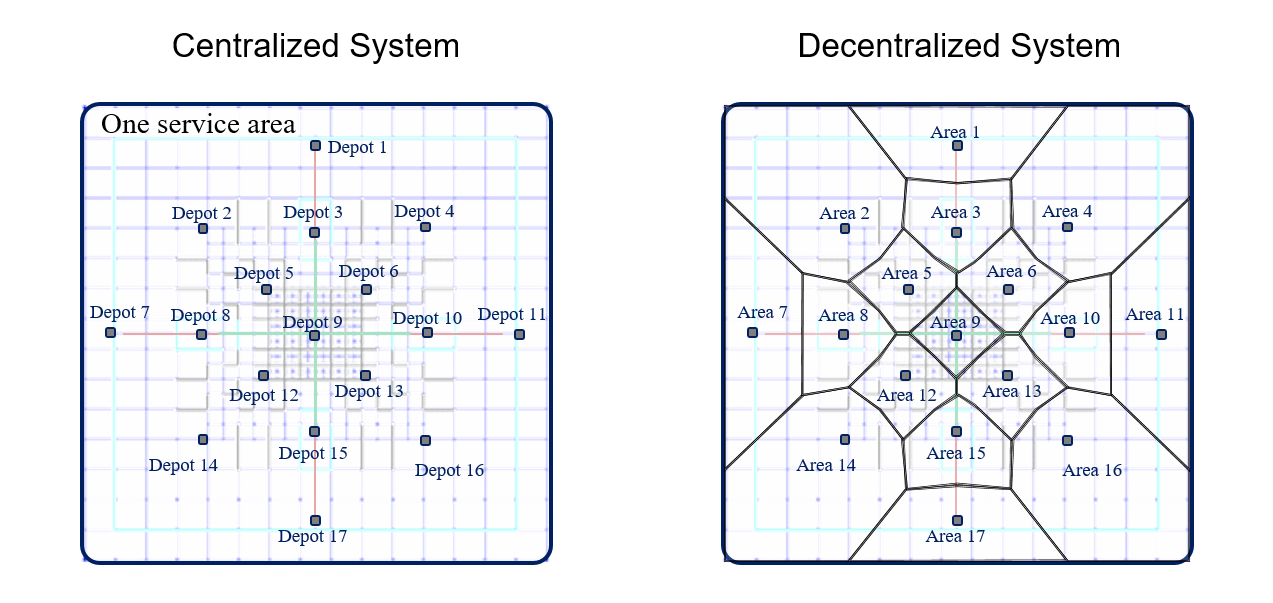}
    \caption{Comparison of centralized and decentralized systems for service area management. The centralized system (left) operates as a single service area with depots distributed across the network, while the decentralized system (right) divides the network into multiple service areas, each managed independently with assigned depots.}
    \label{Fig:Method_Central_Decentral}
\end{figure}

To improve the reliability of RH service operations, various strategies for rebalancing vacant RH vehicles have been explored, primarily through centralized and decentralized systems. In a centralized system, a single main operator manages the entire network with the objective of maximizing overall RH service performance. In contrast, a decentralized system involves multiple local operators, each focused on optimizing the performance of their respective areas independently. The concept for the strategies is shown in Fig.~\ref{Fig:Method_Central_Decentral}.

\subsubsection{Centralized Strategies}
Centralized RH rebalancing strategies have been studied mainly from a driver perspective. For example, \citet{castillo2017surge} and \citet{lu2018surge} suggested dynamic surge pricing strategies that led drivers to decide their relocating area in the next assignment to maximize their profits. They showed improvements in service quality by reducing user waiting times and increasing driver revenue. However, this approach had limitations that certain low-demand areas were not covered due to the indirectness of encouraging drivers to move, which consequently led users in the areas to change their travel options to a more reliable transportation mode. To overcome this issue, many studies have explored the driver decision-making processes when searching for the next passengers. The proposed strategies include minimizing the relocating distance \citep{kim2020vehicle} and the waiting time \citep{meshkani2023centralized}, and maximizing the potential revenue of a driver \citep{afeche2023ride} and the probability of finding a passenger \citep{ao2024control}.

Afterwards, studies have focused on optimizing an agent's sequential decision-making process for a long-term period, and Markov Decision Process (MDP) has gained its popularity \citep{puterman2014markov}. MDP has been used as a stochastic optimization method in which a vacant vehicle takes an action in a given state using its optimal policy for a long-term expected reward. To enhance the RH service performance, MDP has been associated with other methods such as the rolling horizon method \citep{zhou2018optimizing} and queuing theory \citep{yu2019markov}. Furthermore, \citet{shou2020optimal} proposed an e-hailing platform that does not necessarily guide drivers to high-demand areas, but introduces them to low-supply areas with high matching probability. Although these studies showed significant improvements in service quality, there were limitations in that they used predetermined parameters and did not fully capture the traffic dynamics between different transportation modes in road networks. As network traffic can have significant impacts on dynamic RH service operations, it is essential to reflect real-time traffic conditions.

As a consequence, centralized approaches have also been studied considering traffic dynamics (\citealp{goel2016optimal}; \citealp{wang2016pickup}; \citealp{ordonez2017dynamic}). In particular, interest in a trip-based Macroscopic Fundamental Diagram (MFD) model has arisen because of its ability to reproduce the evolution over time of mean traffic conditions for a complete transportation network, using it as a global behavioral curve (\citealp{lamotte2016morning}; \citealp{mariotte2017macroscopic}; \citealp{leclercq2017dynamic}; \citealp{mariotte2019flow}; \citealp{ameli2020cross}). In addition to applying the trip-based MFD model to large-scale networks, they also considered both driver and user perspectives to minimize the loss of ride requests. For example, \citet{alisoltani2021can} introduced a mathematical model that minimized both driver operating costs and user waiting times. \citet{ramezani2023dynamic} proposed a vehicle-passenger matching method using an MFD model. \citet{valadkhani2023dynamic} suggested a non-linear model predictive control framework that relocates idle vehicles in locations with a higher possibility of faster pick-ups.

In the findings of these above-mentioned studies, the improvements in overall performance of RH service operations in the transportation network that capture driver and user perspectives, network congestion, and local supply-demand imbalances can be highlighted. However, they did not account for multi-modality in the network, where users might switch or transfer to other transportation modes if not served. Additionally, scalability remains a significant challenge to address. 

\subsubsection{Decentralized Strategies}
In parallel, decentralized RH rebalancing strategies have been steadily studied, leveraging their ability to scale to large-scale scenarios that involve extensive fleet sizes, a substantial number of ride requests, and complex road networks. In these approaches, individual vehicle entities maximize their objectives based on distributed methods in order to convert a main fleet rebalancing problem into several subproblems that are computationally less expensive. For instance, \citet{guo2018probabilistic} utilized a probabilistic decision approach to optimize each vehicle path to balance supply-demand while minimizing expected cost. \citet{ayala2018spatio} utilized a gravitational algorithm that allows each vehicle to find its relocation area to maximize its own self-interest using incomplete historical information. \citet{hu2019effective} proposed a region partitioning method for vehicle relocation considering historical information of spatial and temporal ride requests, showing high-performance results even when the fleet size reaches 10,000. \citet{buchin2019sampling} utilized a discrete probability approach with each vehicle operating independently, based on historical data for pick-up and drop-off. They suggested that not sharing information with other vehicles may make system operations difficult to manage local supply-demand imbalances. 

Similar to centralized strategies, traffic dynamics has also been considered in decentralized studies. \citet{seppecher2023decentralised} suggested an auction-based fleet rebalancing strategy considering uncertain future ride requests. Each driver independently aimed to find the most profitable local area utilizing a two-sided matching algorithm between drivers and local areas. The trip-based MFD model was applied to reproduce the evolution over time of mean traffic speed by region. As extended work, the operational efficiency of the system was further investigated in terms of uncertainty in the expected number of ride requests \citep{seppecher2023decentralised3} and variations in the number of relocation depots \citep{seppecher2024auction}. Moreover, \citet{beojone2023dynamic} developed a dynamic multi-region MFD model that addresses both congestion alleviation and promotes more sustainable mobility for RH services. \citet{beojone2023relocation} proposed an extended Markov chain method to reduce the number of unserved ride requests using a regional speed-MFD model. \citet{zhu2024coverage} suggested a coverage control algorithm to proactively relocate idle vehicles to minimize user waiting times under traffic dynamics. These aforementioned decentralized strategies demonstrated improvements in system performance and computational power. Moreover, they further considered more details that might occur in real-world situations, such as uncertainty in demand prediction, variation in the number of relocation depots, and noise in the number of background traffic.

Although multiple transportation modes have been included in the network, similar to centralized strategies, multi-modality that allows users to transfer between public transportation modes has not been considered. Providing users with options to choose multiple modes within a single trip not only reflects real-time traffic conditions as local supply-demand balances of different transportation modes constantly vary, but also enables evaluating the validity of the number of active RH vehicles operating in the given situation. Furthermore, the studies have tested strategies only in regular scenarios and have not considered the resilience of the systems during disruptions. To improve the utility of RH services, it is necessary to evaluate the resilience-enhancing role of RH service as a supplementary mode in the interconnected multi-modal transportation systems.

\subsubsection{Reinforcement Learning-based Strategies}
In recent years, machine learning and deep learning methods have gained popularity as powerful techniques for decision-making in RH rebalancing strategies, leveraging effective data-driven techniques to enhance RH system performance \citep{wen2024survey}. Among these, reinforcement learning (RL) has been applied as an effective method in optimizing vehicle agents' decision-making policies. In a dynamic environment where agents do not obtain global observations, RL algorithms such as Q-learning, Deep Q-Network, and temporal difference learning have been employed to derive the optimal policy \citep{ernst2024introduction}. 

RL-based rebalancing strategies have steadily been explored with centralized systems. \citet{shi2019operating} proposed a RL-based centralized decision making framework to convert RH vehicle fleet dispatch problem into a linear assignment problem to minimize user waiting time and operating cost. \citet{liang2021integrated} suggested an efficient RL-based algorithm with a centralized programming module to improve the average profit and reduce matching delay and user waiting time. The RL-based centralized systems showed their strengths in reducing computational burden and improving robustness of the model under system constraints.

However, when training multiple agents with fewer system constraints, these RL-based algorithms are typically limited to small-scale environments due to limitations in computing power. As the environment has become more complex with a group of agents, multiple agents interact with both the shared environment and other agents. As a consequence, multi-agent reinforcement learning (MARL) has been explored to establish decentralized systems, while focusing on agent communication and coordination (\citealp{bucsoniu2010multi}; \citealp{lowe2017multi}; \citealp{yang2018mean}; \citealp{yu2023decentralized}). It is known as a powerful technique due to its success in handling complex and high-dimensional tasks such as playing the game of Go \citep{silver2017mastering}, Poker \citep{brown2019superhuman}, and StarCraft II \citep{vinyals2019grandmaster}.

MARL-based rebalancing strategies have primarily been studied under the distributed nature of the peer-to-peer RH rebalancing problem within decentralized systems. \citet{li2019efficient} proposed a decentralized solution that agents share local information to alleviate the supply-demand gap and increase accumulated driver income. \citet{yang2020multiagent} suggested a proactive model that predicts demand in regions and then relocates vehicles in advance to meet future needs. \citet{shou2020reward} proposed a mean-field MARL-based strategy using a bilevel optimization problem with the upper level as the reward design and the lower level as a multi-agent system. The results demonstrated improvements in order response rate; thus, it can mitigate traffic congestion. \citet{chen2024competitive} integrated a MARL method with a mesoscopic simulation model to address the competitive spatial-temporal pricing problem.

While most MARL-based RH rebalancing studies have enhanced the operational performance of RH systems, they mainly focused on the movements of RH vehicles, neglecting the interconnections with other transportation modes. To effectively integrate RH services into multi-modal transportation systems, it is essential to consider both traffic dynamics across modes and user travel behaviors. Moreover, these studies often assumed that relocated vehicles could be immediately matched with users upon completing their relocation, with relocations occurring to nearby areas within a single time step. In reality, however, vehicles require time to complete relocations and must remain available to fulfill ride requests during the relocation process.

\subsection{Objective and Contributions}
From the above discussion, three key research gaps can be summarized as follows.
First, most previous studies have not considered the aspect of multi-modality for users. Although multiple transportation modes interacted in the road networks, users were unable to switch modes once they departed on their trips. Limited attention has been given to scenarios where users transfer between different modes or adjust their travel plans in response to unexpected changes. To address this, capturing the traffic dynamics of vehicles and users in large-scale networks is essential to accurately reflect real-time traffic conditions.
Second, the impacts of RH rebalancing strategies on transportation systems have been studied in regular scenarios. As the worldwide RH industry grows, further studies regarding resilience with the integration of multi-modal transportation systems need to be assessed across various disruption scenarios.
Third, most studies have applied historical pick-up and drop-off data. However, uncertainty in spatial-temporal future ride requests must be considered to effectively address local supply-demand imbalances with proactive rebalancing strategies. In particular, disruptions provide further discrepancies between prediction and reality. Deep learning approaches also face limitations in achieving optimal evaluation with demand uncertainty, as they are primarily trained on historically known data. 

To bridge these research gaps, the resilience of RH rebalancing strategies in the multi-modal transportation systems during supply disruptions is assessed in this study. A MARL-based RH rebalancing framework, augmented with deep learning capabilities, has been proposed and developed to minimize matching time. RH vehicles, acting as agents, utilize a centralized training and decentralized execution framework, effectively enhancing system performance while reducing computational burden. Moreover, the RH systems are considered to operate under uncertainties, including demand prediction noise and delay in response to disruptions. Using a 900 $\text{km}^2$ toy urban network, a traffic simulation model is implemented using a trip-based MFD model that includes multiple transportation modes such as trains, metros, buses, RH vehicles, and background traffic. Users are set to opt for various transportation options that best fit their preferences and can also change their choices at any time if the situation changes from what they have expected. With disruption scenarios on train lines, the system performance is evaluated with resilience metrics using average travel time and transfer waiting time.

The main contributions of this paper are as follows:
\begin{itemize}
    \item \textbf{Integration of RH Services into Multi-modal Transportation Systems with Demand Uncertainties}: RH mobility services are incorporated into multi-modal transportation systems where users transfer transportation modes based on their preferences. For RH system operations, uncertainties in demand prediction are considered to evaluate system performance during disruptions. The uncertainties consist of variations in demand prediction noise and delay in response to disruptions. This joint design of the multi-modality and the RH system operations contributes to evaluating impacts of RH rebalancing strategies on a more realistic transportation network.
    
    \item \textbf{Assessment of Resilience from User and Operator Perspectives}: Resilience is assessed through key performance indicators, including user waiting time, resilience metrics, total travel time, and distance, from both user and operator perspectives. By monitoring all vehicle and user trajectories, the study evaluates how different RH rebalancing strategies contribute to enhancing the resilience of multi-modal systems.
    
    \item \textbf{Design of an Efficient MARL-based RH Rebalancing Strategy}: A MARL-based strategy of multi-agent deep deterministic policy gradient (MADDPG) is proposed to capture the stochastic supply-demand dynamics in large-scale networks. This MADDPG method with a centralized training and decentralized execution framework optimizes vehicle cooperative behaviors in a more stable way. To minimize the duration of vehicle vacant time, a higher positive reward is granted for quicker matches. Conversely, a greater negative penalty is applied as vacant time extends. Moreover, relocating vehicles are set to be matched with users whenever ride-requests occur nearby, exposing them to a higher uncertainty during service.
\end{itemize}

The remainder of this paper is organized as follows. Section~\ref{section2} presents the principles of multi-modal transportation networks and resilience assessment. Section~\ref{section3} details the components of the RH systems and explains the proposed MARL-based RH rebalancing strategy. Following this, in Section~\ref{section4}, the experimental results section, shows the computational results from the simulated study and provides comparisons to benchmark models. Finally, Section~\ref{section5} concludes the paper and highlights future research directions.

\section{Resilience of Multi-modal Transportation Systems} \label{section2} 
This section outlines the principles of multi-modal transportation systems and resilience assessment. Considering that road networks include both vehicles and users, Section~\ref{subsection 2.1} details the operations of transportation systems, user travel behaviors, and their interactions. Section~\ref{subsection 2.2} introduces scenarios of supply disruptions on train systems and the metrics used for resilience assessment.

\subsection{Multi-modal Transportation Network}\label{subsection 2.1}
\subsubsection{Transportation System Operations}
\begin{figure}[h]
    \centering
    \includegraphics[width=0.9\textwidth]{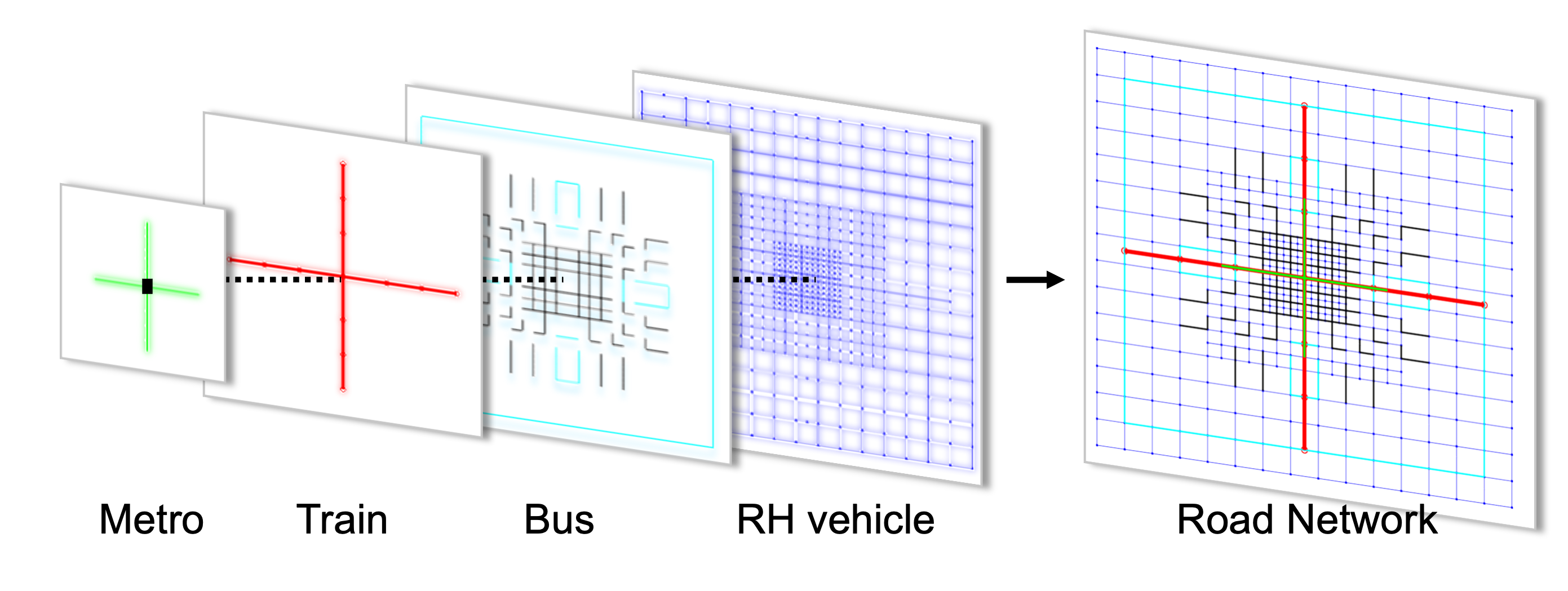}
    \caption{Illustration of individual layers in a multi-modal transportation network, including metro, train, bus, and RH vehicles, and their integration into a comprehensive road network.}
    \label{Fig:Method_Layers}
\end{figure}

For the multi-modal transportation system operations, a road network with multiple transportation mode layers is designed as a basic structure, as shown in Fig.~\ref{Fig:Method_Layers}. Transportation modes include metros, trains, buses, RH vehicles, and background traffic. The layers are connected by nodes, enabling users to transfer for multi-modal trips. 

Mathematical concepts of graphical nodes and edges are used to formulate the road network. The location of the station or the intersection is set as a node, and the link between the two nodes is set as an edge. Mathematical symbols are used to represent the information of the nodes and edges as follows. Let G($\mathcal{P}$,$\mathcal{E}$) be an undirected graph of the network. The set of nodes $\mathcal{P}$ consists of metro stations $\mathcal{P}_M$, train stations $\mathcal{P}_T$, bus stops $\mathcal{P}_B$, and vehicle road intersections $\mathcal{P}_V$. The set of edges $\mathcal{E}$ consists of metro edges $\mathcal{E}_M$, train edges $\mathcal{E}_T$, bus edges $\mathcal{E}_B$, vehicle edges $\mathcal{E}_V$, and connected edges $\mathcal{E}_C$ (connecting nodes of different transportation modes). Here, the edges for transportation modes are obtained from metro line $\mathcal{L}_M$, train line $\mathcal{L}_T$, bus line $\mathcal{L}_B$, and vehicle network $\mathcal{L}_V$. Thus, the edge is expressed as a set of nodes $\mathcal{E}^{i,j}=\{({p}_i, {p}_j) | {p}_i, {p}_j \in \mathcal{P}, {p}_i\neq{p}_j\}$. Metro fleet $\mathcal{F}_M$, train fleet $\mathcal{F}_T$, bus fleet $\mathcal{F}_B$ are running on respective lines. The passenger capacities of the metro, train, bus, and RH vehicle are set as $Q_M$, $Q_T$, $Q_B$, and $Q_V$, respectively. Let $\tau^{i,j}$ be travel time or walking time from node $\mathit{p}_i$ to $\mathit{p}_j$. 

Upon the network design, each transportation mode is operated under the following conditions. The metro and train systems are not affected by traffic conditions, and their average speeds are set at 13 m/s and 30 m/s, respectively. They run the operations based on their given public transportation lines. If the remaining time to reach the next station is calculated to be within a unit time, metros or trains stop by the station in the next time step and stay for another time step, allowing users to have time to ride on or get off the modes. When they reach the final stations of the lines, metros and trains are set to run in opposite directions at predefined times.

In contrast, in bus and car systems, a trip-based MFD model is employed as they interact on the road. In general, travel time of vehicle $i$ of type $m$, denoted as $T_i^m$, is associated with the travel distance $l_i^m$, departure time $t_i$, and the speed $v_m^r (t)$ at time $t$, and it can be mathematically expressed as follows:
\begin{equation}
    l_i^m = \int_{t_i}^{t_i + T_i^m} v_m^r (t) \, dt.  \label{eq:travel_distance} \\
\end{equation}

In the MFD framework, it is assumed that the speed of a transportation mode $v_m^r (t)$ is the same for all users in the same region sharing the same mode at the same time \citep{balzer2022modal}. This speed corresponds to the multi-modal MFD curves, which usually depend on the accumulation number of all transportation modes on the road in the region. The speed formula for buses, private cars, and RH vehicles is defined as follows:
\begin{equation}
    v_m^r (t) = V_m^r \left( \{ n_{m'}^r (t), \forall m' \in \mathcal{M} \} \right) \\
\end{equation}
\noindent where, \( n_{m'}^r (t) \) denotes the number of vehicles of type \( m' \) in service area \( r \) at time \( t \). \( \mathcal{M} \) represents all vehicle types considered in the simulation. \( \{ n_{m'}^r (t), \forall m' \in \mathcal{M} \} \rightarrow \{ V_m^r, \forall m' \in \mathcal{M} \} \) describes the multi-modal MFD speed function, which captures the congestion dynamics and interactions between modes in service area \( r \) at time \( t \). The speed is assumed to be always greater than zero to avoid a complete gridlock. Therefore, the relationship for cars and buses can be formulated as follows using Equation~\ref{eq:travel_distance}:
\begin{flalign}
    & l_i^{car} = \int_{t_i}^{t_i + T_i^{car}} V_m^r \left( \{ n_{m'}^r (t), \forall m' \in \mathcal{M} \} \right) \, dt.  \label{eq:car_travel_distance} \\
    & l_i^{bus} = \int_{t_i}^{t_i + T_i^{bus}} V_m^r \left( \{ n_{m'}^r (t), \forall m' \in \mathcal{M} \} \right) \, dt.  \label{eq:bus_travel_distance}
\end{flalign}

The bus system is operated in the same manner as the metro and train systems. With given lines, buses stay at bus stops for one time step for boarding or alighting users. However, in the RH system, RH vehicles can freely cross the road without any lines. They are operated by transitioning vehicle states between idle, relocating, pick-up, and serving. Further description of the RH services is detailed in Section~\ref{section3}.

\subsubsection{User Travel Behavior} \label{User travel behavior}
\begin{figure}[h]
    \centering
    \includegraphics[width=1.0\textwidth]{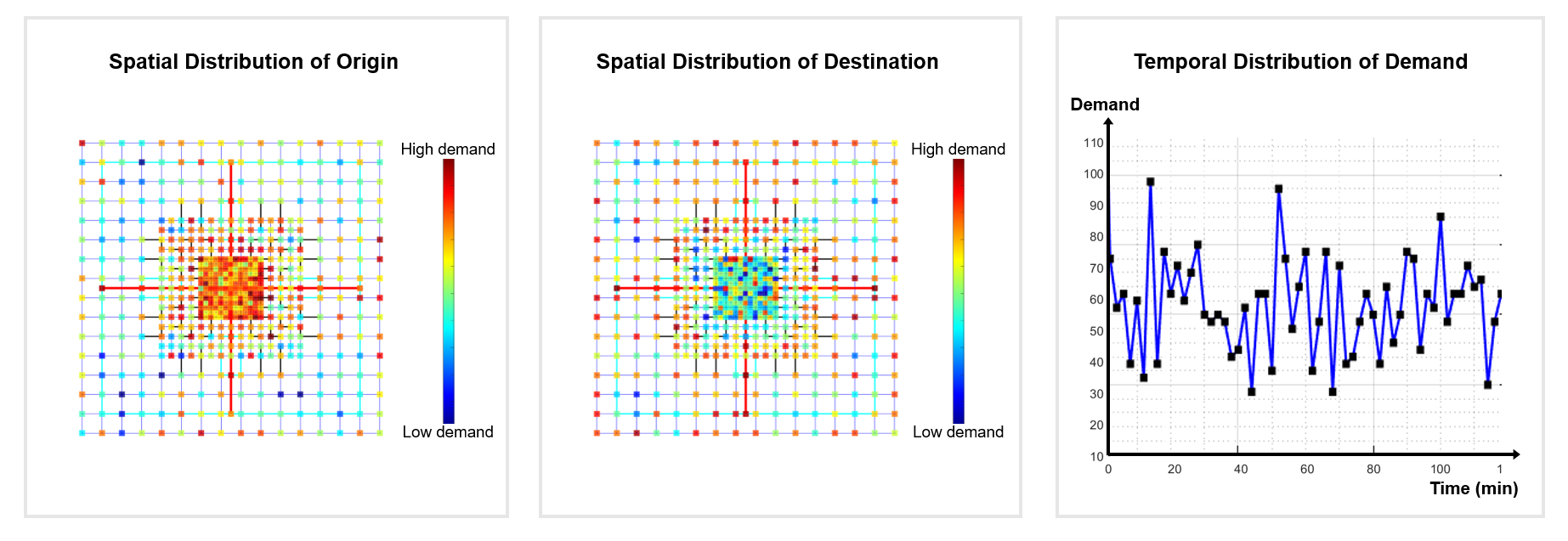}
    \caption{Spatial and temporal distributions of user origin-destination (OD) demand. The left panel shows the spatial distribution of origins, the middle panel illustrates the spatial distribution of destinations, and the right panel depicts the temporal variation of demand over time.}
    \label{Fig:Method_ODdistribution}
\end{figure}

User origin-destination (OD) distribution in this study follows evening commute patterns, with most users traveling from urban to suburban areas. Fig.~\ref{Fig:Method_ODdistribution} shows the spatial and temporal distributions of the user trips. For spatial distributions, urban areas have a high probability as origins and a low probability as destinations, while suburban areas have the opposite probability pattern. For temporal distributions, the time intervals between consecutive user departures follow an exponential distribution with a mean $\lambda$ that applies the Poisson distribution. With this generation, each user is assigned a trip information of when they enter, where they enter, and where they arrive at the end.

To start trips, users select their travel paths in two steps. Firstly, at most five candidate shortest paths are searched using the Dijkstra method. Considering total travel time as a cost, users also consider possible external delays such as walking time, transfer waiting time, and service rate. The candidate paths have different combinations of transportation modes, but parts of them can be the same. Secondly, one candidate is selected as a user path by employing a logit model that describes the user's travel behavior. The utility function, its component, and the probability function of path choice are formulated as follows:
\begin{flalign}
    & U_i = V_i + \epsilon_i   \label{eq:U_i} \\
    & V_i = \sum_{q_i} \frac{d_t (q_i)}{v_t} + \sum_{p_i} \frac{d_m (p_i)}{v_m} + \sum_{l_i} \frac{d_b (l_i)}{v_b (l_i)} + \sum_{s_i} \frac{d_r (s_i)}{v_r (s_i)}   \label{eq:V_i} \\
    & P_i = \frac{e^{-\theta \times U_i}}{\sum_i e^{-\theta \times U_i}}    \label{eq:P_i}
\end{flalign}

\noindent where, Equation~\ref{eq:U_i} represents the random utility of path \( i \), \( U_i \), consisting of the deterministic component \( V_i \) and the unobservable error \( \epsilon_i \). Here, it is assumed that the error term follows a Gumbel distribution. The deterministic component \( V_i \) in Equation~\ref{eq:V_i} is calculated as the total sum of travel times for each transportation mode. Travel times for regions \( q \), \( p \), \( l \), and \( s \) are determined by multiplying the travel speeds \( v_t \), \( v_m \), \( v_b (l_i) \), and \( v_r (s_i) \) by the travel lengths \( d_t (q_i) \), \( d_m (p_i) \), \( d_b (l_i) \), and \( d_r (s_i) \) for trains, metros, buses, and cars, respectively. \( P_i \) in Equation~\ref{eq:P_i} represents the probability that a user will choose the path \( i \).

\subsubsection{Traffic Dynamics of Vehicles and Users}
\begin{figure}[h]
    \centering
    \includegraphics[width=0.95\textwidth]{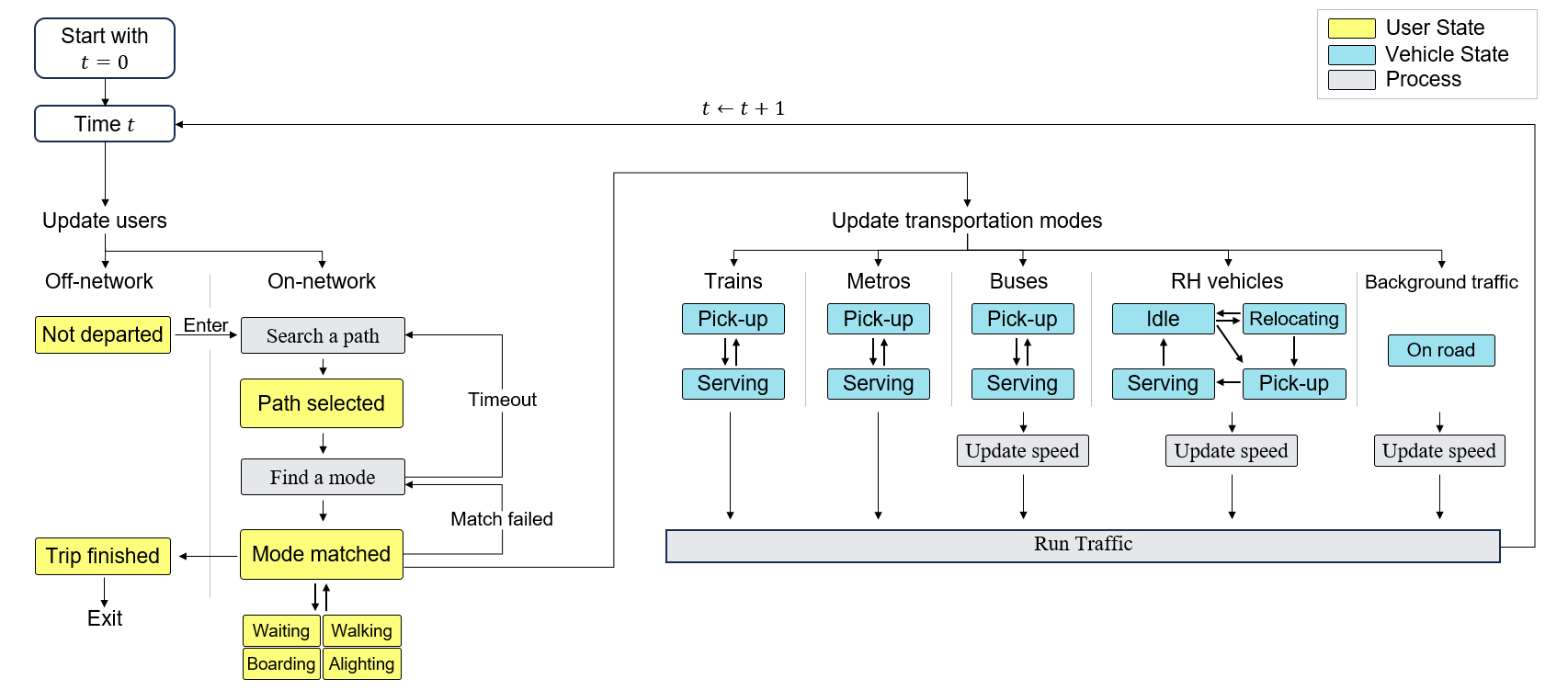}
    \caption{Flowchart of the simulation procedure illustrating traffic dynamics.}
    \label{Fig:Method_TrafficDynamics}
\end{figure}

Upon these principles of vehicle and user behavior, the traffic dynamics for their interactions is designed to simulate a more practical transportation network. A framework of the proposed dynamic simulation procedure is illustrated in Fig.~\ref{Fig:Method_TrafficDynamics}, starting with time $t=0$.

The simulation proceeds as follows. In a given time step $t$, a user pool is updated by distinguishing off-network and on-network users. The off-network users are those who did not yet depart or already completed trips. The on-network users are those who are going to depart in the next time step or are already on a trip. As they enter the network, they search for paths using the path-finding process. After selecting paths, they look for the corresponding transportation modes to be matched. If not matched for a certain amount of time duration, users search for other paths. If matched, they become on-trip users who transition between waiting for a vehicle, walking to a station, boarding a mode, and alighting from a mode. However, if perturbations occur due to unexpected events such as long waiting times or system breakdowns, the matching fails, and these users return to find other available modes to continue their trips. If the matching times out, they search for other paths.

After updating user information, movements of transportation modes are progressed in the same time step $t$. For public transportation, every vehicle of metros, trains, and buses is linked with matched users and updates the information of pick-up and drop-off at each station. For the RH system, each vehicle transitions its state between idle, relocating, pick-up, and serving. The idle state represents a state where the vehicle is vacant and is cruising on the road. The relocating state indicates the state that a vehicle is vacant and relocating to a depot for the next match. These vacant vehicles can be matched with users when nearby ride requests occur. When a vehicle is matched with a user, the state is transitioned to the pick-up state and goes to the user. After pick-up, the vehicle state becomes serving until it drops off the user. Then, the state returns back to idle, and this process repeats. For background traffic, cars on the road are responsible for traffic conditions. After renewing vehicle states, the speeds of the vehicles on the road are updated using the trip-based MFD model. However, metros and trains maintain speeds of 13 m/s and 30 m/s, respectively. Finally, the simulation runs traffic and repeats with a time update of $t+1$.

\subsection{Resilience of Transportation Systems}\label{subsection 2.2}
\subsubsection{Disruption Scenarios}
\begin{figure}[h]
    \centering
    \includegraphics[width=1.0\textwidth]{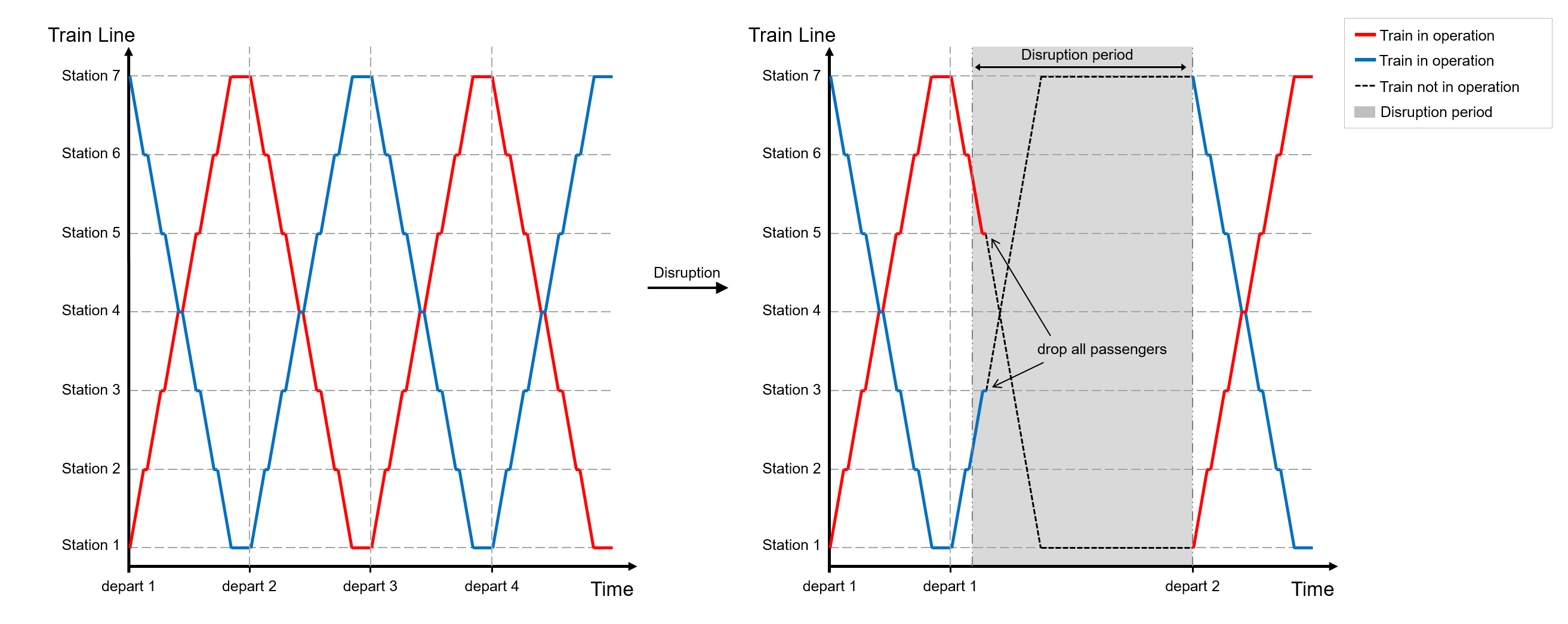}
    \caption{Illustration of an example for train system operations under normal conditions and during a disruption period.}
    \label{Fig:Method_TrainDisruption}
\end{figure}

An illustrative example of a supply disruption of train systems is shown in Fig.~\ref{Fig:Method_TrainDisruption}. In this study, the northbound and southbound train lines are assumed to be connected and trains run in the opposite direction at predefined times after they reach final stations. During a disruption period, all train stations on the affected lines are closed, while the others still operate. Users who ride on affected trains are guided to get off at the following stations regardless of their destination. After dropping off the users, the trains go to the final stations to fix mechanical issues, without allowing any users to board. Furthermore, the network G($\mathcal{P}$,$\mathcal{E}$) transforms into $G'(\mathcal{P},\mathcal{E}')$ representing changes in the availability of the train system. When a disruption ends, the train lines recover their original operations and the trains depart from the starting stations, with the network transforming its state to G($\mathcal{P}$,$\mathcal{E}$).

During disruptions, it is assumed that users do not notice the event in advance. They are affected as their current travel circumstances face the challenge. Stranded users include those who ride on the affected trains and those who just arrived at the affected train stations for the next boarding. For those who planned to board the affected trains but are currently traveling with other modes, they finally realize the disruption as they arrive at the affected stations and the disruption is still ongoing. Following the simulation process in Fig.~\ref{Fig:Method_TrainDisruption}, the stranded users find alternative travel options to finish their trips, which sequentially increases the burden of other transportation modes.

\subsubsection{Resilience Assessment}
The interrelationships among transportation modes are critical components for resilience. During disruptions, users suffer from longer travel times and service delays. With the integration of multi-modal transportation systems, RH service operations can provide a promising way to reduce the burden of other transportation modes and improve resilience. To evaluate the impacts of the RH service, user travel time and transfer waiting time are considered key indicators. The system performance can be assessed by measuring how much travel time is affected by transfer waiting time over time. The mathematical function for the system performance is described as follows:
\begin{equation}
F(t) = \bar{F} - \sum_{\omega \in Q} \sum_{p \in P_w} \Gamma_p (t) \frac{n(P_w)}{n(Q)} 
\end{equation}

\noindent where,  \( \bar{F} \) represents the baseline average travel time during nominal scenarios. \( \Gamma_p (t) \) denotes the user waiting time after being matched with an RH vehicle for path \( p \in P_w \) within the OD pair \( \omega \in Q \). Here, a weighted average of waiting times is obtained that considers variations in OD demand. Then the system performance \( F(t) \) is measured by subtracting the average waiting time from the average travel time.

\begin{figure}[h]
    \centering
    \includegraphics[width=1.0\textwidth]{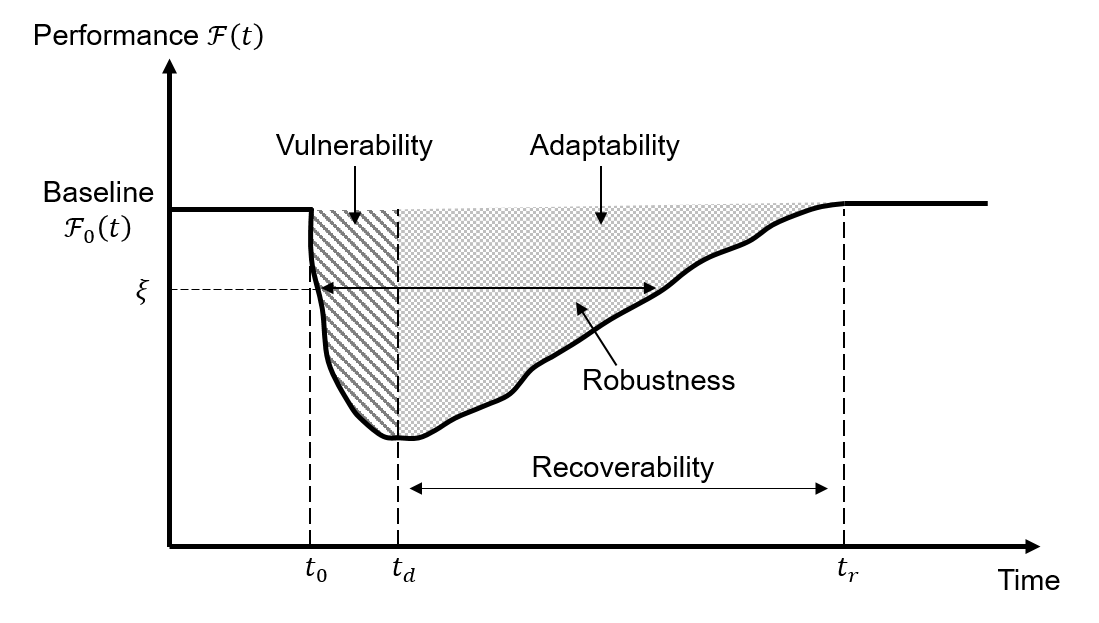}
    \caption{Resilience measurement framework illustrating system performance over time, highlighting vulnerability, robustness, adaptability, and recoverability (adapted from \cite{zhang2024analysis}).}
    \label{Fig:Method_SystemPerformance}
\end{figure}

For a multifaceted assessment of performance, four resilience indicators are defined: vulnerability, adaptability, robustness, and recoverability. The mathematical equations and an illustrative example graph are shown in Equations \ref{eq:R_1}–\ref{eq:R} and Fig.~\ref{Fig:Method_SystemPerformance}, respectively:
\begin{flalign}
    & R_1 = 1 - \frac{\int_{t_0}^{t_d} F(t) \, dt}{\int_{t_0}^{t_d} F_0 (t) \, dt} \label{eq:R_1} \\
    & R_2 = 1 - \frac{\int_{t_d}^{t_r} F(t) \, dt}{\int_{t_d}^{t_r} F_0 (t) \, dt} \label{eq:R_2} \\
    & R_3 = 1 - \frac{H(\xi)}{t_r - t_0} \label{eq:R_3} \\
    & R_4 = \frac{t_r - t_d}{t_r - t_0} \label{eq:R_4} \\
    & R = \alpha_1 R_1 + \alpha_2 R_2 + \alpha_3 R_3 + \alpha_4 R_4 \label{eq:R}
\end{flalign}

\noindent where, vulnerability in Equation \ref{eq:R_1} represents the area gap resulting from a reduction in system performance. This is calculated by integrating the system performance over the time interval from the disruption start point \( t_0 \) to the low peak point \( t_d \). Adaptability in Equation \ref{eq:R_2} represents the area gap that occurs during the recovery in performance. It is calculated over the time interval from the low peak point \( t_d \) to the fully recovered point \( t_r \). Robustness in Equation \ref{eq:R_3} represents the ratio of the time duration when the system performs below a specified level. Recoverability in Equation \ref{eq:R_4} represents the ratio of time duration that measures how quickly the system regained its performance. Resilience in Equation \ref{eq:R} is derived from a weighted sum of these four indicators. In this study, all weights are set to 1. In summary, more resilient systems have lower \( R \) indices.

\section{Ride-hailing Services} \label{section3} 
The RH services are operated in two processes: matching RH vehicles with passengers and rebalancing the vehicles. Since this study does not focus on the matching process, a widely used approach is adopted that matches ride-requesting users with the nearest vacant vehicles. This section presents the RH service system in terms of vehicle rebalancing. Section~\ref{subsection 3.1} introduces preliminaries on the entities of the RH system: system operators and drivers. Section~\ref{subsection 3.2} proposes the MARL-based RH rebalancing strategy.

\subsection{Preliminaries}\label{subsection 3.1}
\subsubsection{System Operators}
The primary objective of system operators is to balance supply and demand in service areas. To this end, they predict user ride-request demand in the near future and request drivers to relocate in advance to secure supply.

Demand prediction is one of the fundamental components for operators that determines the number of relocation vehicles in the areas. Discrepancies between prediction and actual demand worsen operational efficiency. This study starts from the basic assumption that the operators know all user departure times and origin points. Despite that, it is hard to predict the exact number of RH ride-requesting users due to the multi-modality in the network. Taking into account the RH market share, the operators estimate that 10\% of the departing users will be using the RH service \citep{sikder2019uses}. They estimate the expected number of users in the next five minutes, which is the average relocating time for RH vehicles to depots. The definition of demand prediction is as follows.

\begin{definition}[\textit{Demand Prediction}] \label{def:demand_prediction}
A demand prediction \(k\) is the predicted number of ride-requesting users in a given area in the next five minutes.
\end{definition}

Operators cannot deploy additional vehicles whenever there is a shortage. Therefore, it is critical for the system to operate and handle local supply-demand imbalances within the limited number of vehicles. Considering the discrepancies between predicted and actual demand, it is inefficient to always relocate the vehicles to the same number as the predicted demand. When demand becomes higher, the probability that all predicted users actually request is exponentially low. To estimate the probability of matching, a survival function is adapted from \cite{seppecher2023decentralised2}, as shown in Equations~\ref{eq:demand_prediction} and~\ref{eq:demand_prediction_error}.
\begin{flalign}
    & P\left(\hat{k}\right) = \frac{1}{2} \left(1 - \frac{2}{\sqrt{\pi}} \int_0^{\frac{\hat{k} - \mu}{\sigma \sqrt{2}}} e^{-u^2} \, du\right) \label{eq:demand_prediction} \\
    & \hat{k} = k + \epsilon(k,p), \quad \epsilon(k,p) \sim N\left(0, \sqrt{\frac{\pi}{2}}pk\right) \label{eq:demand_prediction_error}
\end{flalign}

\noindent where, Equation~\ref{eq:demand_prediction} represents the probability of the survival function \(P\left(\hat{k}\right)\) that the predicted \(\hat{k}\) users all request RH vehicles. It is defined as the complement of the cumulative distribution function of the normal distribution with a mean \(\mu\) and a standard deviation \(\sigma\) of user demand. Equation~\ref{eq:demand_prediction_error} represents the number of predicted users obtained by adding 10\% of departing users in the next five minutes, \(k\), with a noise \(\epsilon(k,p)\) where \(p\) is the average relative noise. This noise, or uncertainty, in the demand prediction values can even worsen the discrepancies.

\begin{figure}[h]
    \centering
    \includegraphics[width=0.5\textwidth]{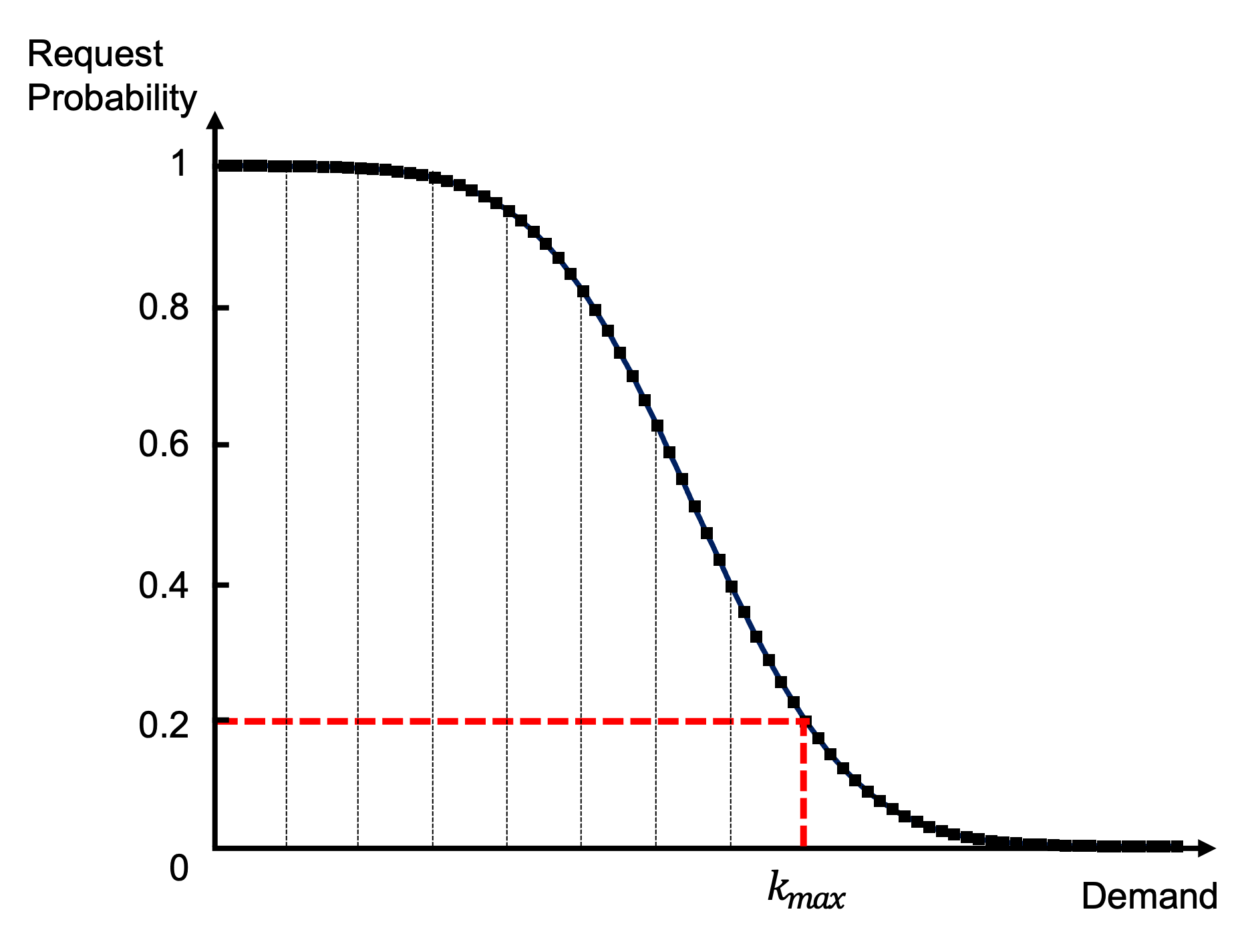}
    \caption{Demand prediction curve illustrating the relationship between request probability and demand, with \(k_{max}\) representing the maximum demand threshold.}
    \label{Fig:Method_DemandPrediction}
\end{figure}

This survival function helps estimate the probability of a completion of the matching for those relocated. If the probability is below a certain threshold level as shown in Fig.~\ref{Fig:Method_DemandPrediction}, it would not be appropriate to publish that amount of relocation offers since the actual number of ride requests might not occur as expected. Operators should prevent a situation where there is local over-supply and under-supply. In this study, the threshold probability is set as 0.2, which is correspondingly paired with the upper bound in the number of relocation offers, \(k_{max}\). In other words, the operators can publish \(k_{max}\) offers at maximum, indicating that the probability of the predicted demand \(\hat{k}\) is lower than 0.2. As a consequence, the number of relocation offers, \(d\), is determined to be \(\hat{k}\) or \(k_{max}\) using Equation~\ref{eq:number_relocation}.
\begin{equation}
d = \underset{x \in \{\hat{k}, k_{max}\}}{\arg\max} P(x) \label{eq:number_relocation}
\end{equation}

However, if disruptions occur, the number of relocation offers may exceed the upper bound to handle a surge in demand in specific areas. Specifically, the estimated number of users stranded on the affected lines is considered for relocation offers.

As operators aim to secure a number of vehicles to handle predicted demand, their utility increases when the number of drivers who decided to relocate is the same as the number of relocation offers. The operator utility function is mathematically expressed as follows.\\

\noindent Operator utility function
\begin{equation}
U_o(r) = |D_r| \label{eq:U_o}
\end{equation}

\noindent subject to,
\begin{flalign}
    & D_r = \{ v \mid Y_v(r) > 0, \ Y_v(r) \leq Y_w(r) \ \forall \, w \neq v, \ Y_v(r) \sim N(\mu_v(r), \sigma_v(r)) \} \label{eq:D_r} \\
    & \mu_v(r) = \frac{d(o_v, r)}{\bar{v}_v} \label{eq:mu_v} \\
    & \sigma_v(r) = \rho \cdot \mu_v(r) \label{eq:sigma_v} \\
    & |D_r| \leq c(r) - i(r) \label{eq:|D_r|1} \\
    & |D_r| + i(r) \geq q(r) \label{eq:|D_r|2}
\end{flalign}

\noindent where, \\  
\hspace*{0.85cm} \( U_o (r) \): operator utility for service area \( r \) \\  
\hspace*{0.85cm} \( D_r \): a set of vehicles to be rebalanced to service area \( r \) \\  
\hspace*{0.85cm} \( Y_v (r) \): a random variable for RH vehicle \( v \) for arrival time to service area \( r \) \\  
\hspace*{0.85cm} \( c(r) \): capacity of depot in service area \( r \) \\  
\hspace*{0.85cm} \( i(r) \): number of vehicles already in service area \( r \) \\  
\hspace*{0.85cm} \( q(r) \): expected ride requests in service area \( r \) \\  
\hspace*{0.85cm} \( \mu_v (r) \): average arrival time to service area \( r \) for RH vehicle \( v \) \\  
\hspace*{0.85cm} \( d(o_v, r) \): shortest distance from RH vehicle \( v \)'s location \( o_v \) to service area \( r \) \\  
\hspace*{0.85cm} \( \bar{v}_v \): average speed of vehicle \( v \) \\  
\hspace*{0.85cm} \( \sigma_v (r) \): standard deviation of average time to arrive at service area \( r \) for RH vehicle \( v \) \\  
\hspace*{0.85cm} \( \rho \): sigma factor (scaling factor for standard deviation) \\

Equation~\ref{eq:U_o} represents the operator utility function for service area \( r \), \( U_o (r) \), as the number of vehicles to be relocated. Equation \ref{eq:D_r} defines a set of vehicles to be relocated, \( D_r \), where arrival times, \( Y_v (r) \), are positive and the closest to a depot of the area \(r\) among the remaining vehicles. The vehicle arrival time of vehicle \( v \) to service area \( r \), \( Y_v (r) \), is modeled to follow a normal distribution with a mean \( \mu_v (r) \) and a standard deviation \( \sigma_v (r) \). Equation \ref{eq:mu_v} calculates the average arrival time \( \mu_v (r) \) by dividing the shortest distance \( d(o_v, r) \) by the average speed \( \bar{v}_v \). Equation \ref{eq:sigma_v} defines a standard deviation \( \sigma_v (r) \) as a proportion of the mean, which is set at 0.1 for this study. Equation \ref{eq:|D_r|1} constrains the number of rebalancing vehicles to not exceed the depot’s vacancy, calculated as the difference between capacity \( c(r) \) and the number of vehicles already in the area \( i(r) \). Equation \ref{eq:|D_r|2} constrains the total number of rebalancing vehicles and those already in the area to be not less than the expected ride requests \( q(r) \).

\subsubsection{Drivers}
As another entity of the RH system, the drivers count on the expected revenue and might refuse the relocation offers if not seen as profitable. Therefore, the driver utility function with the expected profit can be mathematically expressed as follows.\\

\noindent Driver utility function
\begin{equation}
U_v (r) = g_r (k) - c^{\text{rebal}} (o_v, r) \label{eq:U_v}
\end{equation}

\noindent subject to,
\begin{flalign}
    & g_r (k) = \bar{g}_r \cdot max \{P\left(\hat{k}\right), P\left(k_{max}\right)\}\label{eq:g_r} \\
    & \bar{g}_r = \frac{\sum_j d(r_{\text{depot}}, j) \cdot c^{\text{ride}} (r_{\text{depot}}, j)}{\sum_j d(r_{\text{depot}}, j)} \label{eq:bar{g}_r} \\
    & c^{\text{rebal}} (o_v, r) = c^{\text{mileage}} \cdot d(o_v, r) \label{eq:c_rebal}
\end{flalign}

\noindent where,  \\ 
\hspace*{0.85cm} \( U_v (r) \): utility of RH vehicle \( v \) for service area \( r \)  \\ 
\hspace*{0.85cm} \( g_r (k) \): expected earnings in service area \( r \) with minimum \( k \) ride requests  \\ 
\hspace*{0.85cm} \( \bar{g}_r \): average earnings of all combination routes in service area \( r \)  \\ 
\hspace*{0.85cm} \( P\left(\hat{k}\right) \): the probability of ride requests with predicted demand \( \hat{k} \)  \\ 
\hspace*{0.85cm} \( P\left(k_{max}\right) \): the probability of ride requests with the maximum threshold \( k_{max} \)  \\ 
\hspace*{0.85cm} \( X_r \): a random variable of ride request occurrence in area \( r \) \((x \in X)\)  \\ 
\hspace*{0.85cm} \( d(r_{\text{depot}}, j) \): shortest distance from depot \( r_{\text{depot}} \) to location \( j \) in service area \( r \)  \\ 
\hspace*{0.85cm} \( c^{\text{ride}} (r_{\text{depot}}, j) \): riding fare from depot \( r_{\text{depot}} \) to location \( j \) in service area \( r \)  \\ 
\hspace*{0.85cm} \( c^{\text{rebal}} (o_v, r) \): rebalancing cost from vehicle’s location \( o_v \) to service area \( r \)  \\ 
\hspace*{0.85cm} \( c^{\text{mileage}} \): mileage cost per meter  \\ 
\hspace*{0.85cm} \( d(o_v, r) \): shortest distance from RH vehicle \( v \)’s location \( o_v \) to service area \( r \)  \\ 

Equation \ref{eq:U_v} defines the utility function of RH vehicle \( v \) for service area \( r \), \( U_v (r) \), as the profit calculated by the difference between expected earnings \( g_r (k) \) and estimated rebalancing cost \( c^{\text{rebal}} (o_v, r) \). Equation \ref{eq:g_r} expresses that the expected earnings \( g_r (k) \) are obtained by multiplying the average earnings \( \bar{g}_r \) by the probability of ride requests \( P\left(\hat{k}\right) \). Here, \( P\left(\hat{k}\right) \) can be replaced by \( P(X \geq k) \), indicating the probability of exceeding \( k \) ride requests, prompting drivers to accept the operators’ offer to rebalance. Equation \ref{eq:bar{g}_r} calculates the average earnings from all combinations of serving routes in the area. It is obtained considering the ride fare and the shortest distance. Equation \ref{eq:c_rebal} defines the rebalancing cost, which is calculated by multiplying the mileage cost by the shortest distance from the vehicle’s location to the service area.\\

A driver's decision on whether to relocate or not is applied only in the case when the vehicle is vacant. The vehicle that is matched with a user has no other options but to serve. In the real world, however, drivers are sometimes matched with other users for the next turn although they are still serving. For simplicity, the vehicle states are defined as follows.

\begin{definition}[\textit{Vehicle states}] \label{def:vehicle_states}
Vehicle states consist of \textit{Idle}, \textit{Relocating}, \textit{Pick-up}, and \textit{Serving}. Each state has a unique next step, and transitions only when certain conditions are met. Vehicle drivers cannot have a queue of states in advance.
\end{definition}

According to given situations, RH vehicles have different states which transition over time as they operate the service. Four vehicle states are defined in the system: \textit{Idle}, \textit{Relocating}, \textit{Pickup}, and \textit{Serving}. \textit{Idle} represents a vacant vehicle staying at a location without any movements. \textit{Relocating} state indicates a vacant vehicle relocating back to its assigned depot. \textit{Pick-up} denotes a vehicle moving to pick up the matched user. \textit{Serving} represents a vehicle serving the user to the destination.

\begin{figure}[h]
    \centering
    \includegraphics[width=0.8\textwidth]{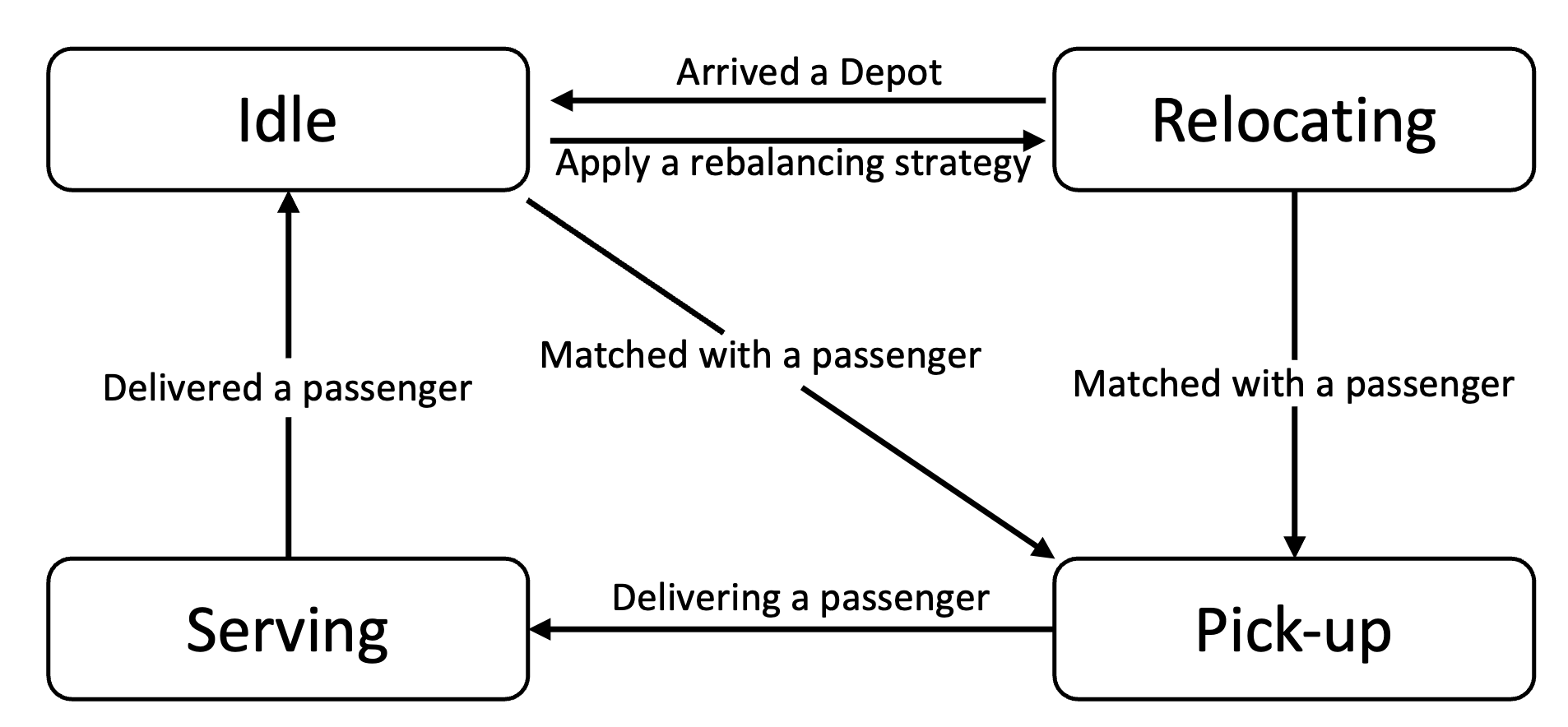}
    \caption{States of RH vehicles and their transitions, including idle, relocating, pick-up, and serving phases. Arrows indicate transitions triggered by events.}
    \label{Fig:Method_RHveh_States}
\end{figure}

Fig.~\ref{Fig:Method_RHveh_States} depicts transitions of RH vehicle states. \textit{Idle} vehicles can maintain the state or transition into \textit{Relocating} if they decide to relocate to depots. \textit{Relocating} vehicles follow shortest paths to reach their assigned depots. When they arrive, they transition states back into \textit{Idle} and wait for the next service. Vehicles in either of \textit{Idle} or \textit{Relocating} are vacant and waiting for the next match with users. If matched, vehicles transition states into \textit{Pick-up} and follow the shortest paths to users' location to pick up. As they pick up users, the state transitions into \textit{Serving}. \textit{Serving} vehicles follow shortest paths to deliver users to their destination. After serving, vehicles become vacant and transition into \textit{Idle}, staying at a location where they dropped users. This process iterates during the entire simulation time.

\subsection{MARL-based Rebalancing Strategy}\label{subsection 3.2}
In this paper, we propose a simulated environment and design a MARL-based rebalancing strategy to address spatio-temporal local supply-demand imbalances under disruptions. The proposed approach enables RH vehicle agents to determine their relocation actions at each time step, specifically when they are \textit{idle}. The following sections detail the key components of this approach, including the design of agents, states, actions, and rewards, the modeling framework, the actor-critic structure, and the optimization process.

\subsubsection{Design of Agent, State, Action, and Reward}
In the simulation environment, the agents observe their states in global and local information. Based on their state, they execute actions and obtain rewards for the actions. The MARL algorithm utilizes the action, state, and reward information to update the decision-making policy. The objective of the policy is to minimize the average vacant time of the agents. A detailed explanation of agent, state, action, and reward is as follows.

\textbf{\textit{Agent}}. Each agent represents a vehicle in a RH system. The agent's primary role is to decide whether to remain idle or relocate to a specific depot to balance supply and demand effectively.

\textbf{\textit{State}}. Each agent observes its environment state only when it is \textit{Idle}. The environment state, represented as \(x = [x_g, x_l]\), consists of global and local information derived from the operating system and the agent's surrounding conditions.

The \textbf{global information} (\(x_g = [g_d, n_d^{RH}, n_d^{PT}, b_d]\)) includes multiple factors that provide a comprehensive view of the system’s supply-demand dynamics:
\begin{itemize}
    \item \textbf{Expected revenue distribution for relocation to depots} (\(g_d = [g_d^1, g_d^2, ..., g_d^n]\)), enabling agents to assess profitability at various locations.
    \item \textbf{Distribution of vehicles already relocated to depots} (\(n_d^{RH} = [n_d^{RH1}, n_d^{RH2}, ..., n_d^{RHn}]\)), helping agents understand current competition at each depot.
    \item \textbf{Availability of alternative transportation modes near depots} (\(n_d^{PT} = [n_d^{PT1}, n_d^{PT2}, ..., n_d^{PTn}]\)), providing additional context on supply alternatives.
    \item \textbf{Supply-demand gap distribution at each depot} (\(b_d = [b_d^1, b_d^2, ..., b_d^n]\)), which highlights areas with potential mismatches that need attention.
\end{itemize}

This global information equips agents with a holistic understanding of system-wide conditions, including supply-demand imbalances, expected revenues based on predicted demand, and uncertainties in user matching. By incorporating data on alternative supply sources, such as nearby ride-hailing vehicles and other transportation modes, agents can make more informed decisions about relocation.

The \textbf{local information} (\(x_l=[n_v]\)) includes the number of \textit{Idle} vehicles within a 500 m distance. Since they are under similar circumstances, this observation provides agents with information to recognize surrounding agents, which prevents an over-supply that can happen if relocated to the same depot. It promotes a cooperative environment with each other \citep{lowe2017multi}.

\textbf{\textit{Action}}. Upon the given state, only \textit{Idle} agents are assigned to take actions for selecting depots for relocation. The set of possible actions for relocation is defined as \( a=[p_0, p_1, p_2, ..., p_k]\), where [\(p_1, p_2, ..., p_k\)] represents the \(k\) possible depots for relocation, and [\(p_0\)] represents the option to cruise without selecting any depot. The action space is set to be continuous to utilize the deterministic gradient framework.

When agents take actions for relocation, their vehicle state transitions from \textit{Idle} to \textit{Relocating}. Once they arrive at the relocation depots with the shortest paths, the vehicle states transition back to \textit{Idle}. If they decide not to relocate to any depots and remain roadside, their vehicle state stays as \textit{Idle}. Since RH vehicles are vacant when in either the \textit{Idle} or \textit{Relocating} states, they can be matched with users when nearby ride requests occur.

\begin{figure}[h]
    \centering
    \includegraphics[width=0.7\textwidth]{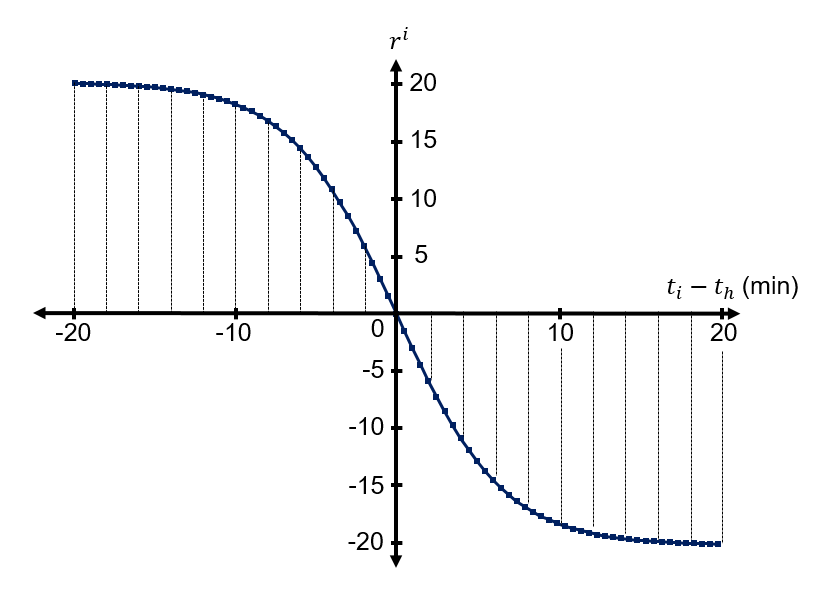}
    \caption{MARL reward function based on vacant time, encouraging quick matching and penalizing delays.}
    \label{Fig:Method_MARL_RewardFunction}
\end{figure}

\textbf{\textit{Reward}}. A reward in MARL is a scalar value that provides the immediate benefit or penalty that an agent receives for taking a particular action in a particular state of the environment. In this study, the reward is designed based on the duration of vacant time following an agent's action. A higher positive reward is granted for quicker matches within a specified time threshold. Conversely, a greater negative penalty is applied as vacant time exceeds the threshold. A reward function for agent \(i, r^i\) is defined in Equation~\ref{eq:reward function} and illustrated in Fig.~\ref{Fig:Method_MARL_RewardFunction}, as follows:
\begin{equation}
r^i = 20 \times \left( 2 \times \frac{1}{1 + e^{0.3 \times (t_i - t_h)}} - 1 \right) \label{eq:reward function}
\end{equation}

\noindent Where, \( t_i\) is the duration that agent \(i\) remains vacant, and \(t_h\) denotes the time threshold. This reward structure for proactive relocation not only aims to minimize users' waiting time but also improves overall system efficiency. For example, quicker matches within the time threshold result in a more negative \(t_i - t_h\), allowing agents to earn a higher positive reward of up to a maximum of 20. However, if matches occur beyond the time threshold, \(t_i - t_h\) becomes positive, leading to a greater negative penalty, down to a minimum of -20.

\subsubsection{Modeling Framework}
\begin{figure}[h]
    \centering
    \includegraphics[width=1.0\textwidth]{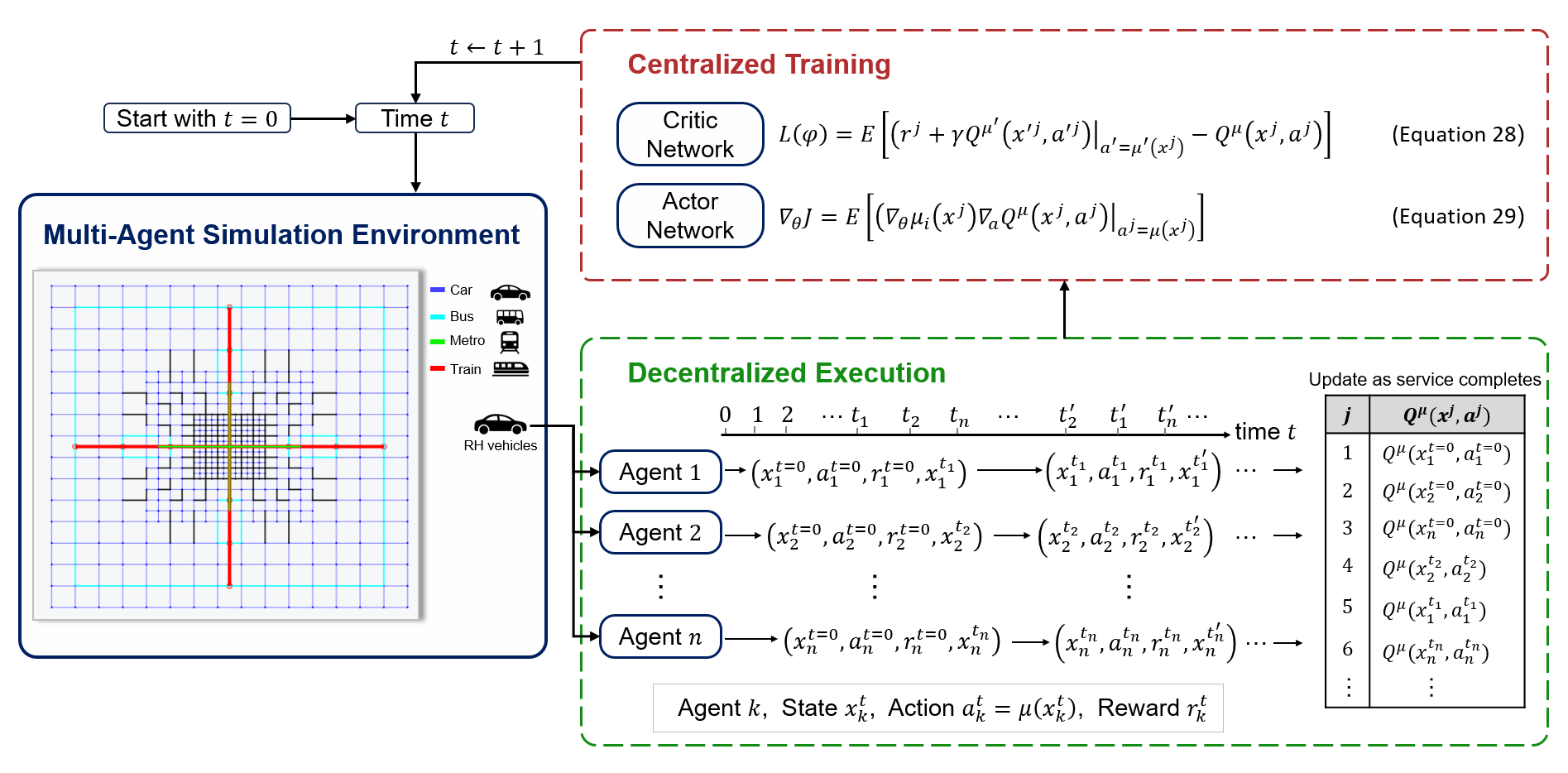}
    \caption{MARL framework with centralized training and decentralized execution (CTDE) for independent agent actions in a simulated environment.}
    \label{Fig:Method_MARL_Process}
\end{figure}

In multi-agent settings, independently learning policy gradient methods such as Q-learning and DQN perform poorly, as each agent's policy changes to overwhelm the others during training, resulting in a non-stationary environment \citep{lowe2017multi}. To address this issue and achieve convergence of the algorithm, the centralized training and decentralized execution (CTDE) framework is adopted as shown in Fig.~\ref{Fig:Method_MARL_Process}. 

Each agent is made up of an actor and a critic. The actor determines the agent's action, which is the output of a decision-making policy with inputting the agent's state. These actions are taken independently under different local states while also sharing global states, which is known as decentralized execution. The main advantages of this approach include reducing increasing scalability and operational efficiency of the system, and adapting to traffic dynamics. Moreover, we propose agents to take actions not necessarily at the same time step, as they are independently serving passengers with different OD pairs under different regional traffic conditions. 

On the other hand, the critic evaluates if the decision-making policy is optimized using a reward function, by inputting the state observed and the action taken by the agent. These critics from all agents are integrated in one dataset for training, which is known as centralized training. By doing so, the optimal policy is applied to each agent to maximize global objectives, or the overall rewards. This method allows agents to consider the potential decisions and impacts of other agents during training, facilitating more coordinated behavior patterns.

This actor-critic structure in the CTDE framework leverages the strengths of both policy and value-based methods. Continuous action spaces paired with decentralized execution and centralized critics provide an environment to develop cooperative behaviors between agents. The decentralized execution ensures scalability and flexibility, while the centralized critic ensures that the agents' actions align toward achieving cooperative goals.

\subsubsection{Actor-critic Structure}
\begin{figure}[h]
    \centering
    \includegraphics[width=1.0\textwidth]{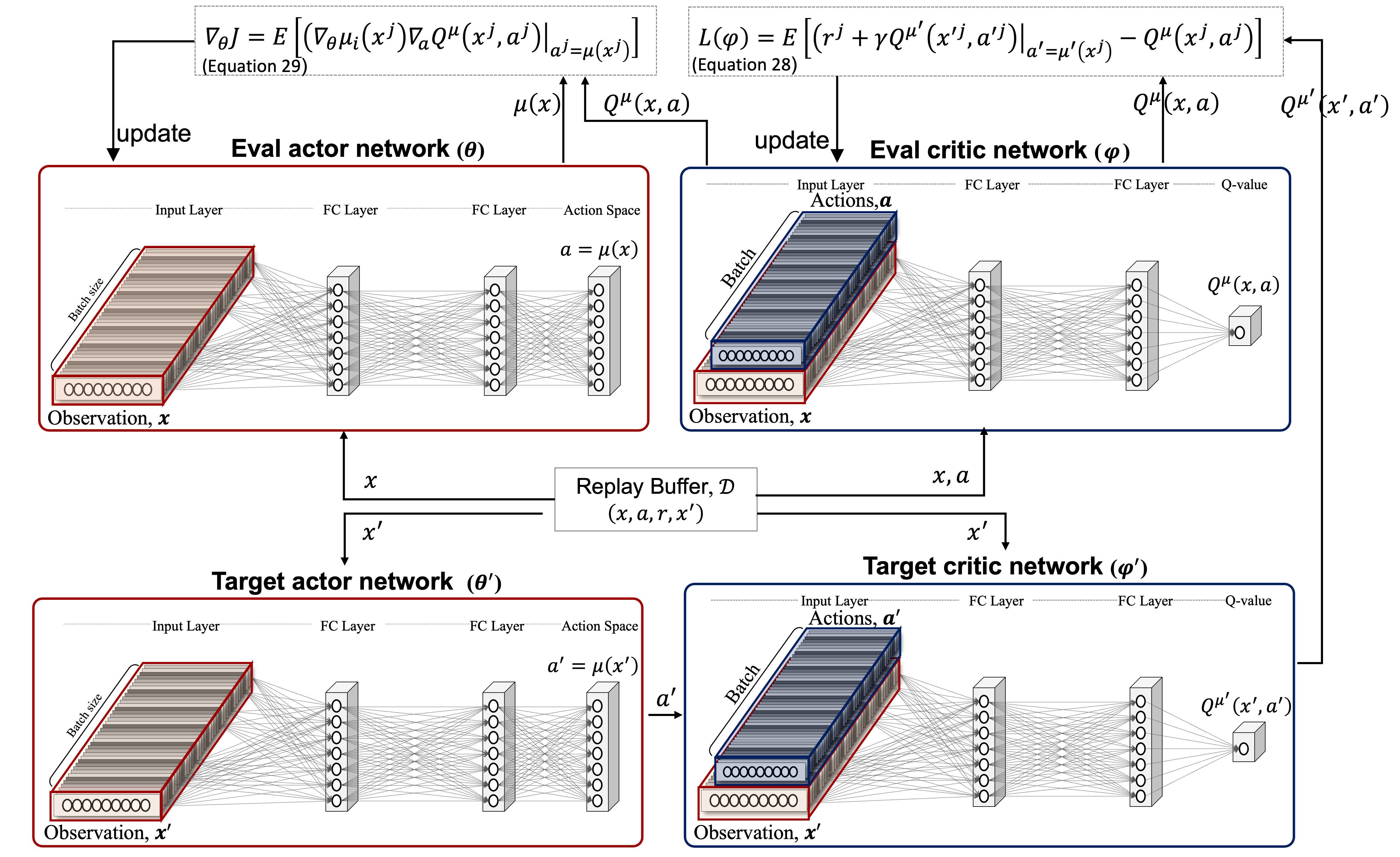}
    \caption{Actor-Critic structure with evaluation and target networks for policy and value updates.}
    \label{Fig:Method_MARL_ActorCritic}
\end{figure}

Using the CTDE framework, all training data are stored in one replay buffer of experiences, enabling agents to share their experiences to optimize the policy. As depicted in Fig.~\ref{Fig:Method_MARL_ActorCritic}, each agent utilizes actor and critic networks in the optimization process. Each agent inputs its environment state to the actor network for decision-making. The critic network evaluates if the agent's action is well-decided using the reward. In this study, the networks are approximated using a deep neural network with two fully connected layers. 

Regarding policy optimization, this study considers a multi-agent extension of MDPs called partially observable Markov games \citep{littman1994markov}. Since this study considers the action space and the policy as continuous, the MADDPG method is adopted. MADDPG utilizes the eval network and the target network in order for the parameters to converge in a more stable way.

More concretely, agent policy parameters are trained by using the eval critic network and the target critic network. The critic networks assess the agent's action \(a\) by reward \(r\) when state \(s\) is given, using a loss function as follows:
\begin{equation}
L(\varphi) = \frac{1}{S} \sum_j \left( r^j + \gamma Q^{\mu'} \left( x'^j, a'^j \right) \big|_{a' = \mu'(s^j)} - Q^{\mu} \left( x^j, a^j \right) \right)^2 \label{eq:Lossfunction}
\end{equation}

\noindent Where, \( L(\varphi)\) is the loss function when evaluating the reward using the eval critic network with a set of parameters \(\varphi\). It is calculated by dividing the total sum of squares of loss by the number of samples \({S}\). \(r^j\) denotes the reward of sample \(j\). $\gamma \in [0,1]$ represents the discount factor. \(Q^{\mu'}\) and \(Q^{\mu}\) are the Q functions of the target critic network and the eval critic network, respectively. \(x'^j\) and \(x^j\) are the next state and the previous state, respectively. \(a'^j\) denotes the action generated by the target actor network \(\mu'\). \(a^j\) is the action of sample \(j\) stored in the historical record.

The eval actor network \(\mu\) outputs action \(a\) with input state \(s\). To optimize decision-making and maximize cumulative rewards, the eval actor network is trained by the sampled policy gradient using the loss function in Equation~\ref{eq:Gradient}:
\begin{equation}
\nabla_{\theta} J \approx \frac{1}{S} \sum_j \nabla_{\theta} \mu_i(x^j) \nabla_{a} Q^{\mu} \left( x^j, a^j \right) \big|_{a^j = \mu(x^j)} \label{eq:Gradient}
\end{equation}

Then, the target networks' parameters are updated by a soft update mechanism using a hyperparameter $\tau \in [0,1]$. This process is implemented to make the parameters softly converge to those of the eval networks under the stationary environment. The parameters of the target networks are updated using those of the eval networks:
\begin{flalign}
    & \varphi' \leftarrow \tau \varphi + (1 - \tau) \varphi'  \label{eq:target_critic_network} \\
    & \theta' \leftarrow \tau \theta + (1 - \tau) \theta'     \label{eq:target_actor_network}
\end{flalign}

\subsubsection{Optimization Process}
In this study, a two-layer 64 unit ReLU multilayer perceptron is used for each network. The Adam optimizer is adopted with a learning rate of 0.01. $\tau$ is set to 0.01 to update the target networks. $\gamma$ is set to 0.99. The size of the replay buffer is 10$^6$. The network parameters are updated every 10 samples added to the replay buffer. Each sample includes a batch size of 500 serving experiences. Five recently updated samples and five random seed samples are selected for each training.

The following Algorithm~\ref{algo:MADDPG} shows how the MADDPG method updates its parameters. Agent executes action for relocation only if it is in the idle state. After it completes serving and becomes idle again, a new state is observed. Since this new state is the next state of the preceding situation, the data including state \( s \), action \( a \), reward \( r \), and the next state \( s' \) are stored in the memory pool \( \mathcal{D} \). Agents' experiences are stored in the memory pool in order as they complete each service, regardless of vehicle IDs. As the memory pool is filled with a sufficient number, samples are randomly extracted for training. In training, using the error measured from the critic eval network (\( Q \) and critic target network (\( Q' \), the parameters of the critic eval network and the actor eval network \( \mu \) are updated. Finally, the target network parameters are updated via a soft update mechanism.

\FloatBarrier
\begin{algorithm}
\caption{Multi-Agent Deep Deterministic Policy Gradient (MADDPG) for Vehicle Rebalancing.} \label{algo:MADDPG}
\begin{algorithmic}[1]
\STATE INPUT: agents, replay buffer \(\mathcal{D}\)
\STATE OUTPUT: target network parameters 
\FOR{timestamp = 1 to \( M \)}
    \FOR{\( t = 1 \) to max-episode-length}
        \FOR{agent \( i = 1 \) to \( N \)}
            \IF{vehicle state is \textit{Idle}}
                \STATE observe state \( s' \) 
                \IF{vehicle has history of previous action \( a \) and previous state \( s \)}
                    \STATE observe reward \( r \) of action \( a \)
                    \STATE Store \( (s, a, r, s') \) in replay buffer \( \mathcal{D} \)
                \ENDIF
                \STATE Select action \( a = \mu_{\theta}(s') + \mathcal{N}_t \) w.r.t. the current policy and exploration
                \STATE Execute action for vehicle rebalancing
                \STATE \( s \leftarrow s' \)
            \ENDIF
        \ENDFOR
        \IF{( n(replay buffer \(\mathcal{D}\)) // \( N \)) > 10}
            \FOR{agent \( i = 1 \) to \( N \)}
                \STATE Sample a random minibatch of \( S \) samples \( (s^j, a^j, r^j, s'^j) \) from \( \mathcal{D} \)
                \STATE Update critic by minimizing the loss using Equation~\ref{eq:Lossfunction}
                \STATE Update actor using the sampled policy gradient using Equation~\ref{eq:Gradient}
            \ENDFOR
            \STATE Update target network parameters for each agent using Equations~\ref{eq:target_critic_network}~and~\ref{eq:target_actor_network}
        \ENDIF
    \ENDFOR
\ENDFOR
\end{algorithmic}
\end{algorithm}
\FloatBarrier

\section{Numerical Results} \label{section4} 
This section assesses how the proposed strategy contributes to the multi-modal transportation systems. Section~\ref{subsection 4.1} represents the design of the experiment for the application. Section~\ref{subsection 4.2} analyzes vehicle dynamics in terms of the number of serving vehicles and the distribution of RH vehicle states over time. Section~\ref{subsection 4.3} assesses system performance with variations in demand prediction noise. Waiting times of different strategies are compared during regular scenarios to assess nominal performance. After that, system performance is further compared during disruption scenarios. Then, the performance is quantified using the resilience index, or R-index. Furthermore, Section~\ref{subsection 4.4} represents the system performance with variations in demand prediction noise and delay in response to disruptions. Lastly, Section~\ref{subsection 4.5} presents how RH rebalancing strategies contribute to the systems by evaluating the relationship between total travel time and total travel distance.

\subsection{Design of Experiment} \label{subsection 4.1}
\begin{figure}[h]
    \centering
    \includegraphics[width=0.8\textwidth]{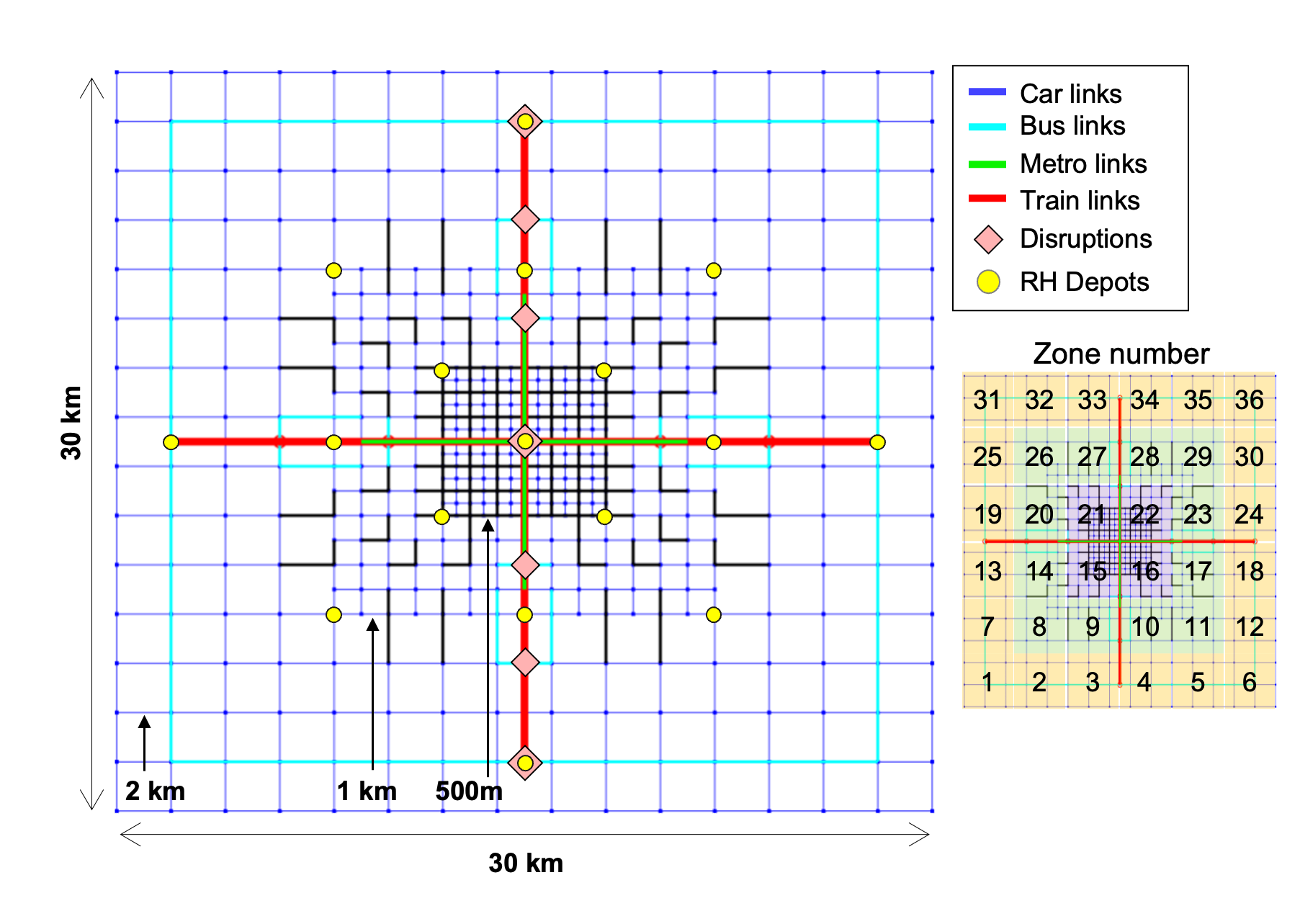}
    \caption{Manhattan grid network illustrating the integration of car, bus, metro, and train links, along with ride-hailing depots (yellow circles) and disruption locations (red diamonds) across a 30 km × 30 km area.}
    \label{Fig:Results_ApplicationSite}
\end{figure}

Vehicle rebalancing strategies are evaluated using a Manhattan grid network, as shown in Fig.~\ref{Fig:Results_ApplicationSite}. It is a widely used road network structure instance for various transportation-related analyses \citep{bayliss2024mobility}. The network covers 900 $\text{km}^2$ with a 30 km-by-30 km grid. It consists of 874 road segments of 500 m, 1 km, and 2 km links, and 576 intersections. The network can be split into 36 zones, each measuring 5 km by 5 km. The urban area includes the 4 innermost purple zones, consisting of 500 m road segments with high accessibility to all transportation modes. The intermediate area surrounds the urban area with 12 green zones, composed of 500 m, 1 km, and 2 km of road segments. It has lower accessibility because the train stations are located far from parts of the area. The suburban area consists of the remaining 20 yellow zones. Most road segments are 2 km links, and only a few zones have access to train stations.

Public transportation includes trains, metros, and buses, with frequencies of 20 minutes, 6 minutes, and 10 minutes, respectively. A total of 600 RH vehicles are operated. Seventeen RH depots are strategically placed for efficient rebalancing: 9 near train stations for joint multi-modal trips, and 8 across the urban area for regular demand. The centralized strategy is managed by a single system operator for all depots. For the decentralized strategy, 17 local operators are managing their assigned areas with their respective depots. Each service area covers nodes that are closest compared to other depots. Lastly, for the MARL-based strategy, the entire network is managed by a single operator, but each vehicle acts based on its given observation.

The simulation time is set for 4 hours with 5,729 users showing evening commute patterns, with most users traveling from urban to suburban areas. The simulation begins with a one-hour warming-up time to keep traffic operations in a stable state. Afterward, a regular scenario is applied for another 30 minutes, with a disruption following for one hour from 1.5 to 2.5 h in timestamp. During the disruption period, seven train stations are closed for northbound and southbound travel. To facilitate agent-based traffic dynamics, Multi-modal Network Modelling and Simulations (MnMS) is applied in this study \citep{MnMs}.

A total of five different rebalancing strategies are tested and compared to evaluate how RH services can improve the resilience of the multi-modal transportation systems. These strategies have been selected because they are the proposed, the state-of-the-art, or the baseline models.
\begin{itemize}
    \item \textbf{Proposed MARL-based strategy}: RH vehicles under the proposed MARL-based strategy communicate and cooperate while independently deciding on relocations based on both global and local observations. This approach shares similarities with the decentralized strategy.
    
    \item \textbf{Centralized strategy}: A single main operator manages the entire network with the objective of maximizing overall RH service performance. The objective function aims to maximize the total sum of operator and driver utilities, adapting the approach of \citet{duan2020centralized}. Details of the algorithm are provided in Appendix~\ref{Appendix:Centralized}.
    
    \item \textbf{Decentralized strategy}: Multiple local operators independently optimize the performance within their respective areas. An auction-based Gale-Shapley matching method, adapted from \citet{seppecher2024auction}, is employed. Details of the algorithm are provided in Appendix~\ref{Appendix:Decentralized}.
    
    \item \textbf{Random rebalancing strategy}: In this strategy, vacant RH vehicles relocate randomly to available depots. As there is no system operator, demand prediction is not utilized.

    \item \textbf{No-rebalancing strategy}: Each vacant RH vehicle remains at the location where it dropped off a passenger, waiting until it is matched with another user. As there is no system operator, demand prediction is not utilized.
    
\end{itemize}

\begin{figure}[h]
    \centering
    \includegraphics[width=0.5\textwidth]{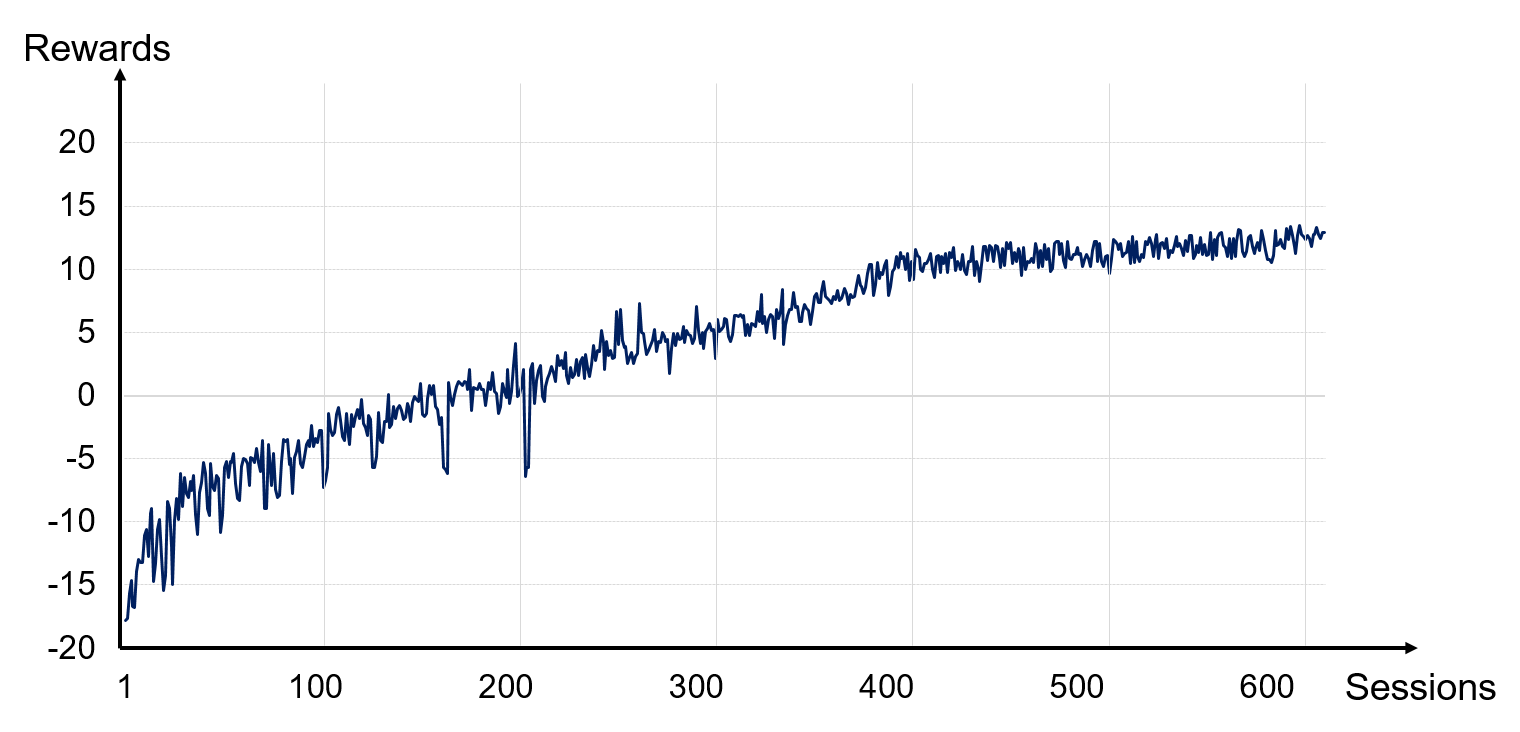}
    \caption{Convergence of the waiting time reward using the MARL-based strategy over training sessions.}
    \label{Fig:Results_MARL_Training}
\end{figure}

For the MARL-based strategy, the network parameters for the decision-making policy were optimized using regular scenarios without disruptions. Fig.~\ref{Fig:Results_MARL_Training} shows the convergence of the average rewards to minimize the waiting time of the RH service with 600 vehicles. The graph illustrates that as the number of training iterations increases, the reward tends to increase and eventually converges at a fixed value. It indicates that the model gradually generates a better policy as the training progresses and ultimately finds the optimal one. After 600 training sessions, the average waiting time is reduced to approximately 4.3 minutes per service. This MARL model is applied to disruption scenarios.

\subsection{Analysis of Vehicle Dynamics Over Time} \label{subsection 4.2}
\subsubsection{Number of Serving Vehicles Over Time}
\begin{figure}[h]
    \centering
    \includegraphics[width=0.8\textwidth]{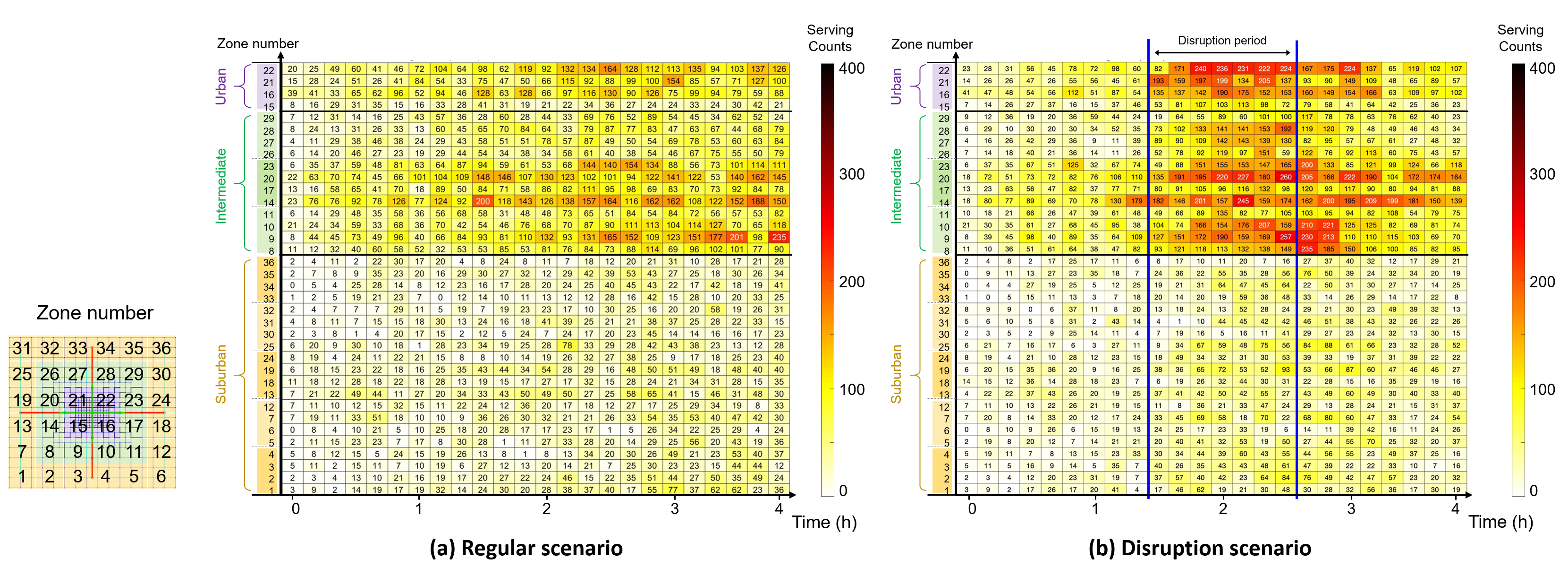}
    \caption{The number of serving vehicles over time in (a) regular and (b) disruption scenarios.}
    \label{Fig:Results_ServingVehicles}
\end{figure}

We first demonstrate the reliability of the simulation model. The results of the proposed strategy in Fig.~\ref{Fig:Results_ServingVehicles} illustrate the number of vehicles serving users across 36 zones over time. In the regular scenario shown in Fig.~\ref{Fig:Results_ServingVehicles}(a), RH mobility services remain inactive until the warm-up period concludes at the 1 h timestamp. Once traffic operations stabilize, the intermediate area appears as the most active service area. This is consistent with the evening commute pattern, where users in the intermediate area have fewer transportation options as they travel from the urban area. High activity levels in zones \#9, \#14, \#20, and \#23 highlight that RH vehicles often serve as last-mile connections, particularly near public transportation stations or stops. Furthermore, the total activity in suburban areas exceeds that in urban areas, despite the suburban area's less-structured road network and limited transportation options. This indicates that ride demand for RH vehicles is higher in the outer areas than in the inner areas during the evening commute. This observation denotes the potential for local supply-demand imbalances when demand predictions are inaccurate, leaving most vacant vehicles roaming urban areas while waiting for the next service requests.

However, under a disruption scenario, the results revealed a completely different pattern. As shown in Fig.~\ref{Fig:Results_ServingVehicles}(b), the number of service vehicles increased across the network at each timestamp. Notably, activities in the urban area increased significantly, exceeding those of the intermediate area. Zone \#22 experienced a peak of 220 serving vehicles for nearly an hour, the highest among all zones. Even after operations returned to normal, zones \#22 and \#16 in the urban area continued to maintain elevated serving numbers for a certain period. These results demonstrate how a disruption on train lines sequentially increases the burden on other modes of transportation, particularly RH services.

In the intermediate area, each zone exhibited patterns similar to those observed in the urban area, but the magnitude of the increase was smaller. Furthermore, two distinct patterns emerged depending on local accessibility. Zones without public transportation stations were affected for shorter periods. For example, servings numbers in zone \#29 gradually increased from 19 to 100 during the disruption and then decreased steadily afterward. In contrast, zones with public transportation stations experienced prolonged impacts. Servings numbers in zone \#14 remained above 150 even after the disruption ended. These different patterns can be attributed to variations in road segment lengths and accessibility to alternative transportation modes. Specifically, one-third of zone \#14 consists of 500-meter links, while only one-ninth of zone \#29 includes such links. Shorter links result in more nodes, thereby increasing the number of users located within the zone. Regarding accessibility, zone \#14 is located near train and metro stations and has more bus routes, including black and blue lines, aligned to the zone. During the disruption, passengers from other modes opted not to transfer to trains. Consequently, it results in a surge in RH vehicle activity in zone \#14. In contrast, zone \#29 lacks connections to trains or metros and has fewer bus routes aligned with it. This indicates that fewer external sources bring RH passengers to zone \#29, resulting in a smaller increase in demand.

In the suburban area, the number of servings remained largely unchanged in most zones, primarily due to the lower number of nodes and limited accessibility. In the case of zone \#34, passengers intending to alight at the final train station could not arrive because of the suspended train operations during the disruptions. These passengers may have chosen alternative routes originating from other zones, bypassing zone \#34. Conversely, those planning to board the train in zone \#34 likely requested RH vehicles as an alternative mode to reach their destinations. The fluctuations in user numbers in the suburban area followed a pattern similar to that of the evening commute pattern, with minor ups and downs in demand.

\subsubsection{Distribution of Vehicle States Over Time}
\begin{figure}[h]
    \centering
    \includegraphics[width=1.0\textwidth]{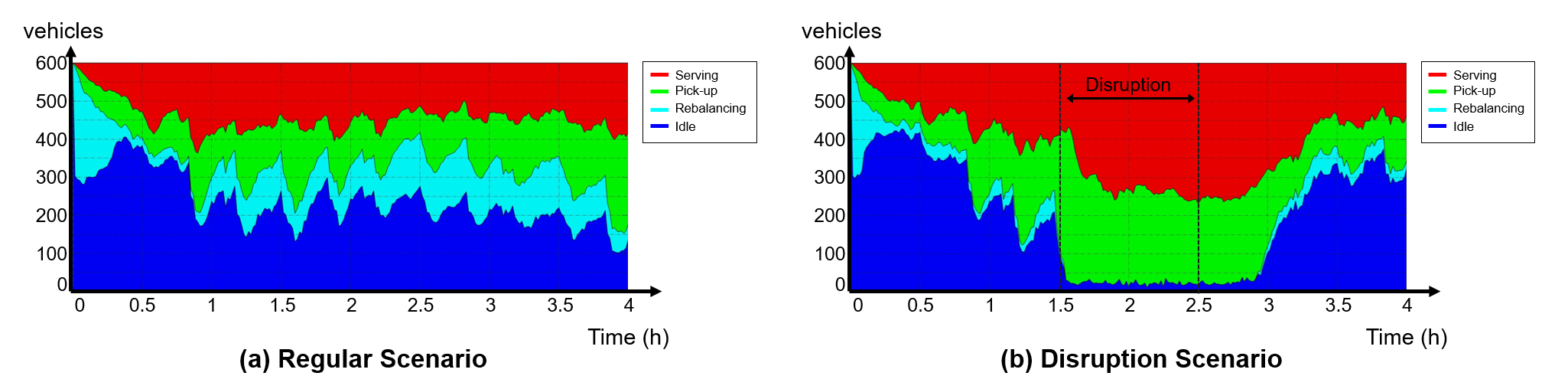}
    \caption{Temporal distribution of vehicle states in (a) regular and (b) disruption scenarios.}
    \label{Fig:Results_VehicleStates}
\end{figure}

We then assess the operational efficiency of the RH service from the operator's perspective. The distributions of vehicle states in both regular and disruption scenarios with the proposed strategy are represented in Fig.~\ref{Fig:Results_VehicleStates}. The x-axis represents time, while the y-axis shows the cumulative number of vehicles. Vehicle activities, including serving, pick-up, rebalancing, and idle, are represented in red, green, cyan, and blue, respectively. The results of the other strategies demonstrated similar patterns.

In a regular scenario, the number of serving and pick-up vehicles increased as they were matched with users during the warm-up period. Simultaneously, the number of idle vehicles increased as RH vehicles were relocated to depots. As the time reached 1 h, the network had attained an equilibrium state, and the number of vehicles in different states began to fluctuate within certain ranges. The number of serving vehicles showed smaller variations compared to the other states. It indicates a consistent user demand for RH vehicles. This consistency aligns with the demand generation process described in Section~\ref{User travel behavior}. The stability in demand generation ensures that the study focuses solely on assessing the impacts of rebalancing strategies on disruptions, without additional influences from demand fluctuations. 

During the disruption, the number of serving vehicles rapidly increased. Nearly all 600 vehicles were engaged in serving or picking up users throughout the disruption period and afterward. This indicates that a significant number of stranded users have made ride requests as an alternative means to complete their trips. The additional demand, combined with the limited availability of other transportation modes, also impacted other users, which led to oversaturation of the RH service. This cascading effect persisted even after the disruption ended. The operational level eventually returned to its original state at the 3.5 h timestamp, one hour after the disruption.

The number of rebalancing vehicles converged to zero during the disruption, indicating that there was no opportunity to perform relocations. This suggests that it is nearly impossible to relocate RH vehicles when the transportation systems has reached its low peak. However, by effectively addressing local supply-demand imbalances, rebalancing strategies can mitigate these challenges. They can help the system performance decline more gradually, improve the low peak, reduce the duration spent at the low peak, and enable faster recovery. 

\subsection{Effects of Demand Prediction Noise on System Performance} \label{subsection 4.3}
Discrepancies between actual demand and prediction can worsen the operational efficiency of transportation systems. Thus, the system performance was evaluated with variations in demand prediction noise by increasing it to \(p=20\%\). The following subsections present comparison results on user waiting times during regular scenarios, system performance during disruptions, and the resilience index.

\subsubsection{User Waiting Times During Regular Scenarios}
\begin{figure}[h]
    \centering
    \includegraphics[width=0.9\textwidth]{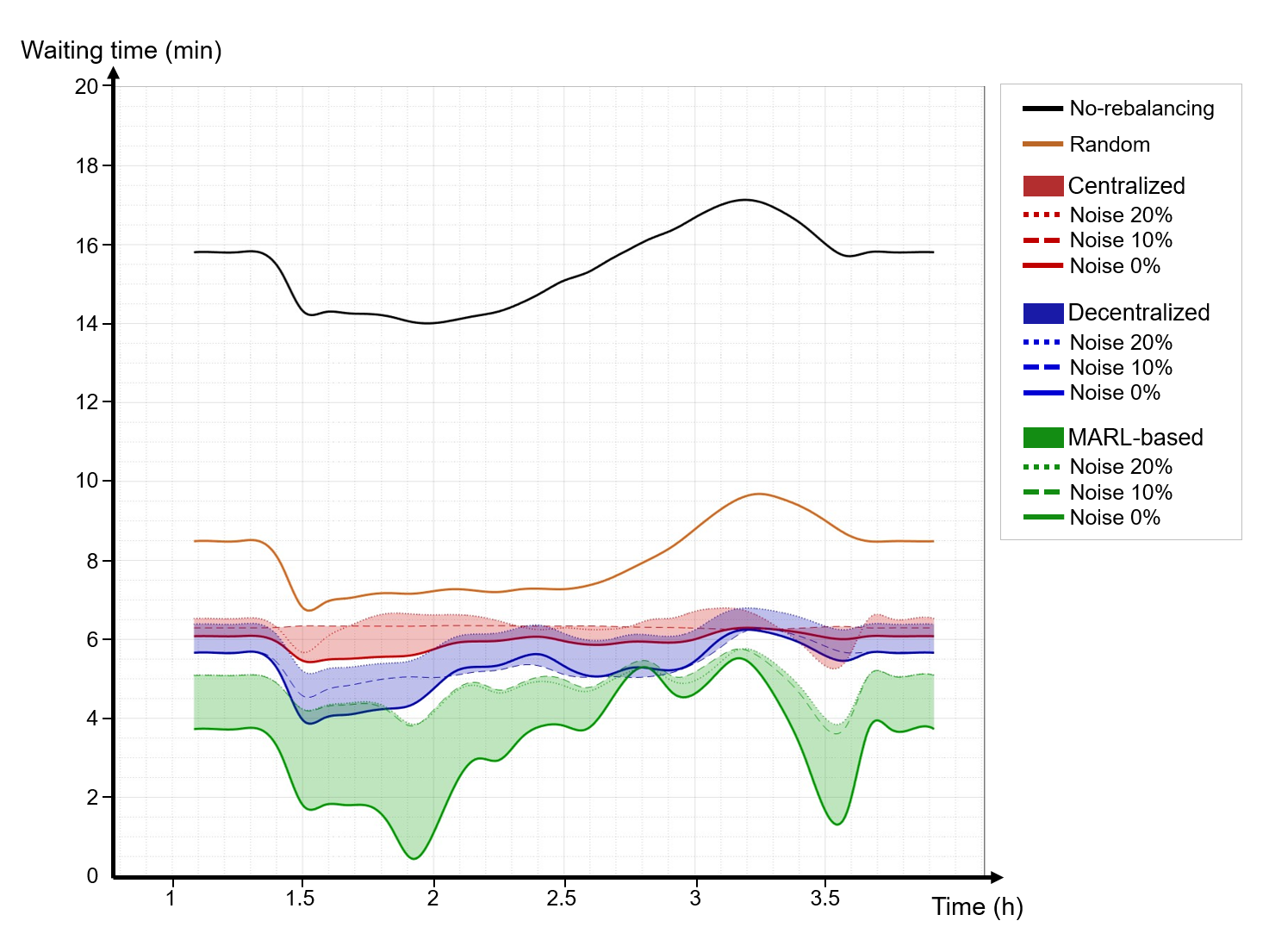}
    \caption{User waiting times with demand prediction noise during a regular scenario.}
    \label{Fig:Results_WaitTime_RegularScenario}
\end{figure}

Before assessing resilience, the system performance of each strategy under regular scenarios is evaluated using the average user waiting time, as shown in Fig.~\ref{Fig:Results_WaitTime_RegularScenario}. The waiting times for each strategy exhibit varying ranges while maintaining respective stable state without disruptions throughout the simulation. For the baseline without the rebalancing strategy, the average waiting times were estimated at 15.4 min, which is the highest among the strategies. The values fluctuated within a range of around 14 to 17 min. However, applying rebalancing strategies is shown to significantly improve operational efficiency. The random rebalancing strategy performed with an average waiting time of 8.0 min, ranging from around 6.5 to 10 min. These two strategies show solid lines, since the prediction noise is not considered.

The remaining strategies showed increases in waiting times as noise levels increased. The centralized strategy showed the most robust performance, with an average waiting time of 6.2 minutes, ranging from approximately 5.5 to 6.5 minutes. The decentralized strategy performed slightly better, with an average waiting time of 5.6 minutes and a range of approximately 4 to 6.5 minutes. This improved performance is attributed to local system operators managing their respective areas, thereby diminishing the impact of noise, which only affects local areas. Lastly, the MARL-based strategy achieved the best performance with an average waiting time of 4.3 minutes. However, it also showed the largest range, ranging from approximately 0.5 to 5.5 minutes. This indicates that the training sessions were insufficient and the decision-making policy is less optimized. Nevertheless, it is important to note that extending the training duration might risk overfitting.

\subsubsection{System Performance During Disruptions}
\begin{figure}[h]
    \centering
    \includegraphics[width=0.90\textwidth]{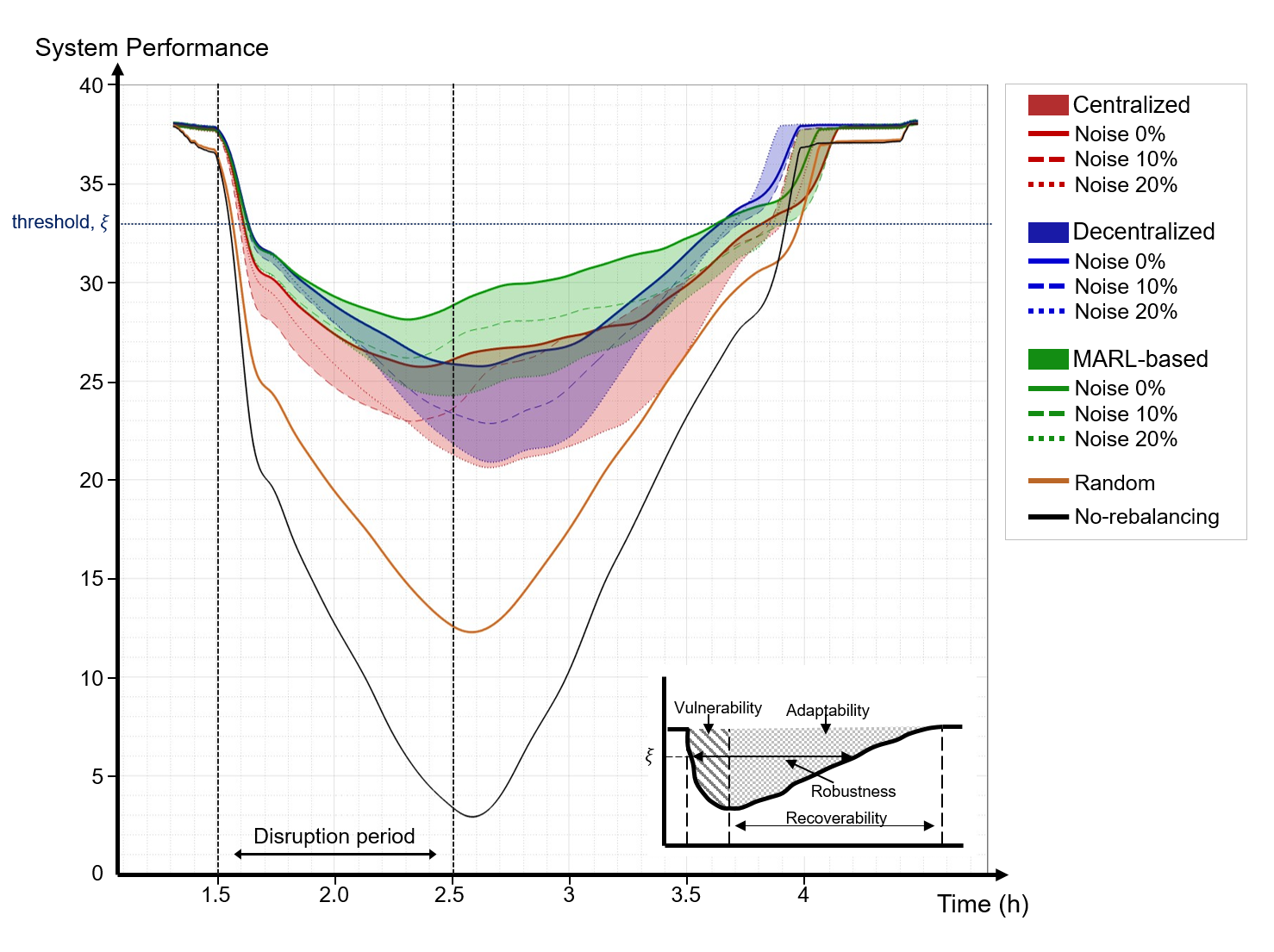}
    \caption{System performance under demand prediction noise during disruptions, comparing resilience across rebalancing strategies.}
    \label{Fig:Results_Performance_DemandNoise}
\end{figure}

Disruption scenarios with demand prediction noise were implemented to evaluate how different strategies contribute to the resilience of the system performance, as presented in Fig.~\ref{Fig:Results_Performance_DemandNoise}. For effective assessment, baseline performances at the start were normalized to the same level. Performance curves show convex patterns that decrease as the disruption occurs and increase after the disruption. In particular, the moments they fully recovered to the baseline were different in time as they applied different strategies. The threshold for robustness, \( \xi \), was set at 85\% of the baseline.

In the case of the no-rebalancing strategy, when the disruption occurred at timestamp 1.5 h, the performance consistently declined at the steepest slope and reached its low peak immediately after the disruption ended. It then recovered with a constant slope, and the slope increased near the threshold. Since no strategies are implemented in the transportation system, the performance declined at an almost constant rate, reached its low peak, and then recovered at another nearly constant rate. Similarly, the random rebalancing strategy showed a similar pattern, but achieved better performance values overall. The application of random rebalancing resulted in significant increases in vulnerability and adaptability, along with small decreases in robustness and recoverability.

The centralized and the decentralized strategies showed both similarities and differences in their performance patterns. Their low peaks tended to decrease as noise levels increased, and both maintained similar low-peak levels under the same noise conditions. However, the centralized strategy tended to reach its low peak during the disruption and recovered slowly to the normal state, displaying a wider convex shape. In contrast, the decentralized strategy reached its low peak after the disruption and recovered more quickly. Consequently, the four resilience indicators for the centralized strategy were measured to be lower than those for the decentralized strategy. Lastly, the MARL-based strategy demonstrated an intermediate performance decline and recovery between the centralized and decentralized strategies. Furthermore, its low peaks consistently showed the highest values among all strategies, regardless of the prediction noise. 

\subsubsection{Resilience Index Assessment}
\begin{figure}[h]
    \centering
    \includegraphics[width=0.90\textwidth]{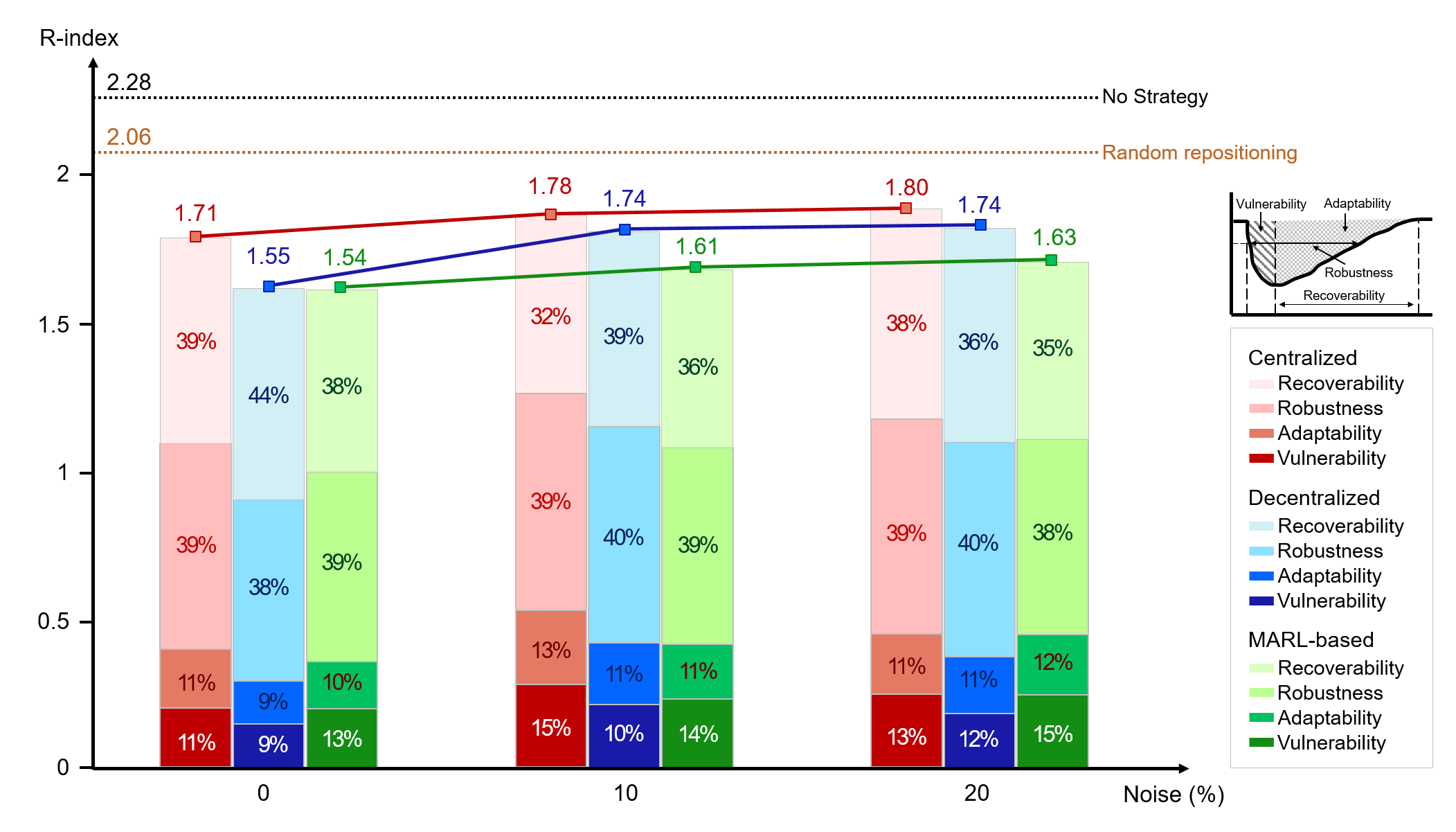}
    \caption{Resilience assessment with rebalancing strategies across different demand prediction noise.}
    \label{Fig:Results_Index_DemandNoise}
\end{figure}

Quantitative evaluation is an effective approach to assessing system performance. Correspondingly, the resilience index of the strategies is compared, as illustrated in Fig.~\ref{Fig:Results_Index_DemandNoise}. The R-index showed an increasing trend as noise levels increased across all strategies, indicating a decline in network resilience. For example, the index for the centralized strategy gradually increased from 1.71 to 1.80 as noise increased from 0\% to 20\%. All four resilience indicators - vulnerability, adaptability, robustness, and recoverability - also increased. These results indicate that greater discrepancies between actual and predicted demand negatively affect system resilience. However, the findings highlight that implementing any rebalancing strategy is more effective than not applying any.

The four resilience indicators presented varying patterns across different levels of noise and strategies. Since lower values indicate high performances, the system showed strong performance in terms of vulnerability and adaptability but weaker performance in robustness and recoverability. For the MARL-based strategy, vulnerability gradually increased from 13\% to 15\% as noise levels increased. Similarly, adaptability increased from 10\% to 12\%. Based on the definitions in Equations \ref{eq:R_1} and \ref{eq:R_2}, these increases imply a worsening of the system's low performance peak. It indicates that stranded users are likely to experience significantly longer wait times for RH vehicles during disruptions. However, robustness remained relatively stable within the range from 39\% to 38\%. This indicates that noise had minimal impact on the relative time it took for system performance to drop and recover concerning the threshold \(\xi\). However, it is important to note that these values are relative to the R-indices and represent the graph patterns; in absolute terms, performance deteriorated further under higher noise levels. Lastly, recoverability decreased from 38\% to 35\%, resulting in a more right-skewed graph. This decline indicates that disruptions lasted longer, and the system's recovery initiation was increasingly delayed.

Based on this resilience assessment, the effects of prediction noise can be summarized as follows. First, an increase in prediction noise worsens overall resilience regardless of the rebalancing strategies. The results of the R-index show that greater improvement in resilience is associated with higher prediction accuracy. Second, an increase in prediction noise increases the ratio of adaptability and vulnerability. The low peak of system performance is deepened, and users stranded in the worst situation may experience longer waiting times for RH vehicles. Third, robustness is maintained, and recoverability decreases as noise increases. This indicates a right-skewed performance graph, where higher noise prolongs suffering until the performance reaches its low peak. Lastly, the centralized strategy showed inconsistent patterns with noise. Given that it measured the highest R-index, the centralized system might not effectively handle disruptions and may even worsen under higher noise levels.

\subsection{Effects of Demand Prediction Noise and Delay in Response on System Performance} \label{subsection 4.4}
System operators manage RH services based on expected demands. Time delays may occur in receiving information about disruptions. If a delay occurs, the cumulative estimated number of stranded users on the affected lines is added to the demand prediction. For example, with a 20-minute delay, the operator relies on regular demand predictions until informed. Once informed, the accumulated surge in users is incorporated into the demand after 20 minutes. According to Equation~\ref{eq:number_relocation}, the number of relocation offers may exceed the upper bound in certain areas to address the surge in demand. System performance was evaluated under combinations of demand prediction noise up to \(p=20\%\) and response delays of up to 30 minutes.
The following subsections present comparison results on system performance during disruptions and the resilience index.

\subsubsection{System Performance During Disruptions}
\begin{figure}[h]
    \centering
    \includegraphics[width=0.9\textwidth]{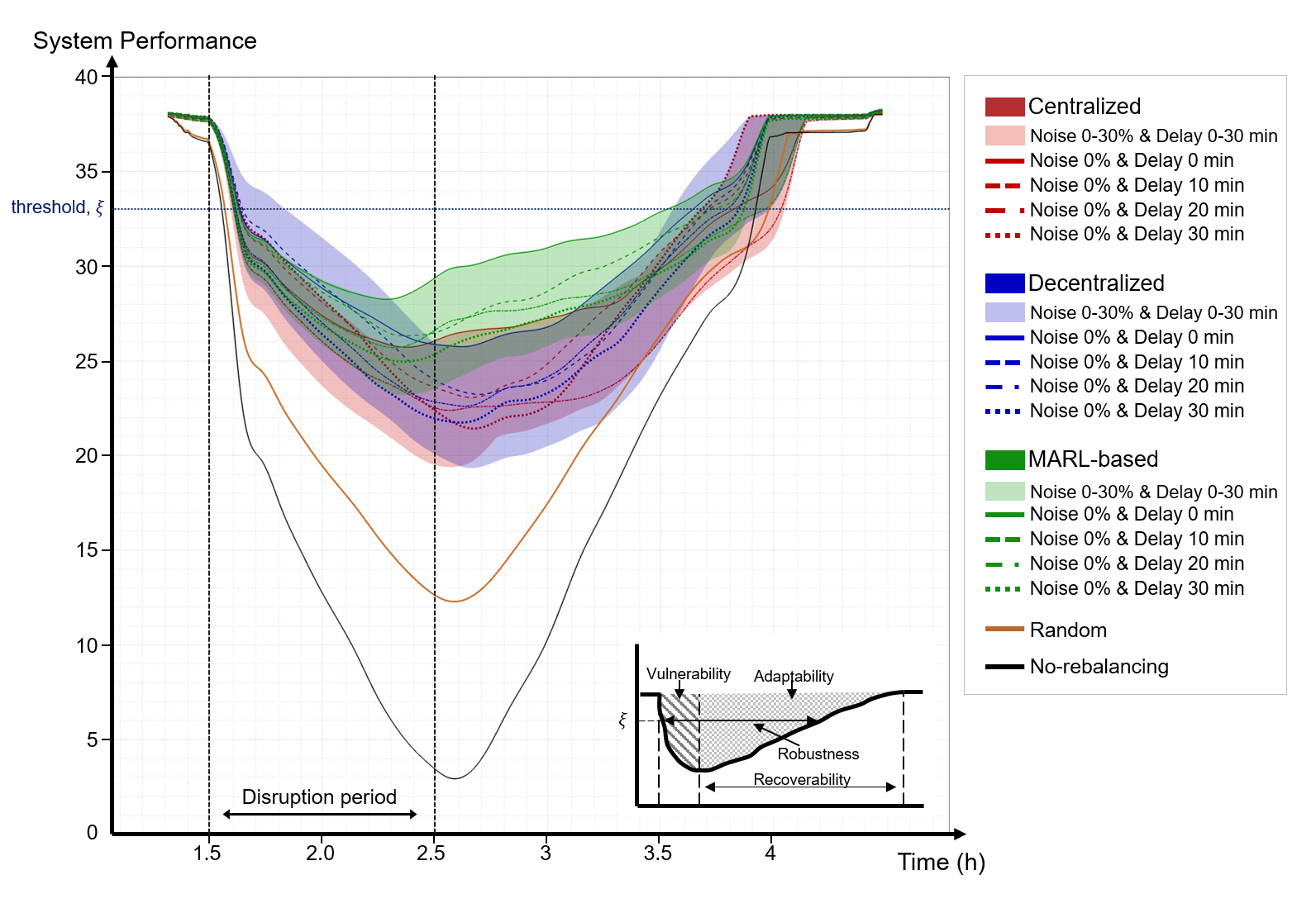}
    \caption{System performance with demand prediction noise and delay in response to disruptions.}
    \label{Fig:Results_Performance_ResponseDelay}
\end{figure}

Disruption scenarios with demand prediction noise and response delays were implemented to assess the resilience of system performance, as shown in Fig.~\ref{Fig:Results_Performance_ResponseDelay}. Uncertainties in demand prediction were not applied in the no-rebalancing and random rebalancing strategies, as these strategies do not involve a system operator.

For the centralized strategy, low performance peaks were observed within the range of approximately 20 to 26. When a disturbance occurred, the variations in time taken to reach the threshold \(\xi\) were similar to those observed under scenarios with only prediction noise. However,during recovery, the variations in time required to recover above the threshold showed a broader range. This suggests that increased uncertainty in demand prediction negatively impacts stranded users, resulting in significantly longer waiting times in the worst cases. Consequently, more users are exposed to disruption for extended periods due to delays in system performance recovery.

The decentralized strategy also experienced a performance decline comparable to the centralized strategy. The time duration with the performance below 25 increased relative to scenarios with only prediction noise. Additionally, the variations in time to drop below and recover above the threshold increased. Nevertheless, the decentralized strategy showed the slowest decline and the fastest recovery times compared to other strategies.

The MARL-based strategy showed the smallest performance decrease, indicating that users experienced the least negative impacts on average waiting time. However, despite its relative advantage, the uncertainties still imposed deeper and broader performance fluctuations, adversely affecting the resilience indicators.

\subsubsection{Resilience Index Assessment}
\begin{figure}[h]
    \centering
    \includegraphics[width=1.0\textwidth]{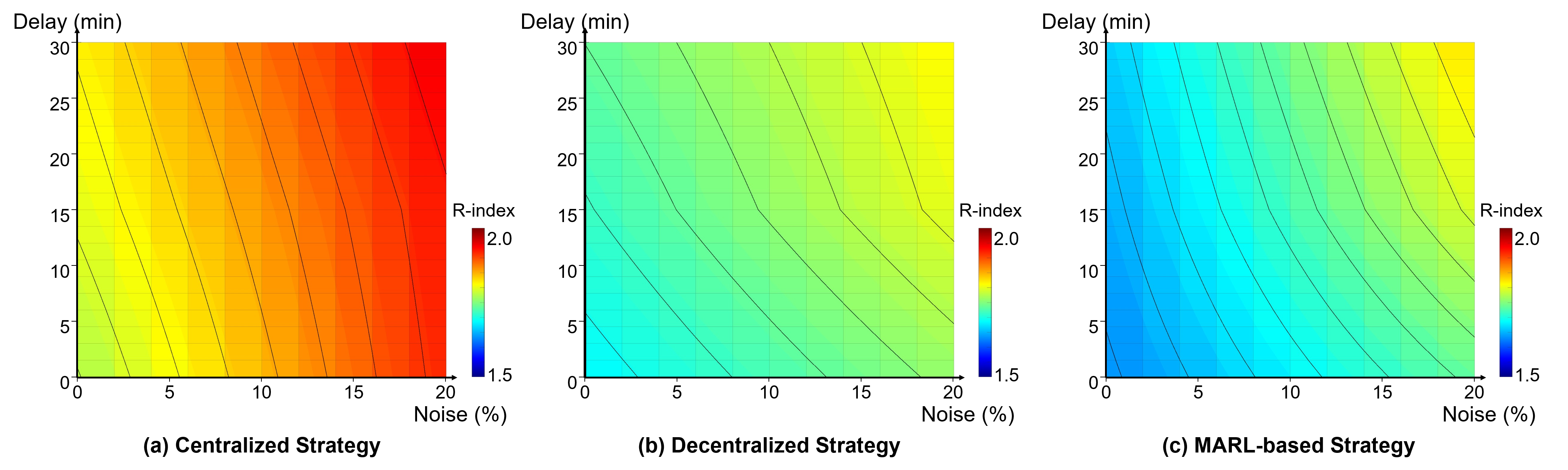}
    \caption{Resilience index with prediction noise and delay in response to disruptions.}
    \label{Fig:Results_Index_ResponseDelay}
\end{figure}

The resilience of the strategies is quantitatively compared, as presented in Fig.~\ref{Fig:Results_Index_ResponseDelay}. The R-index generally increased across all strategies as noise and delay increased. Among the strategies, the centralized system showed the highest R-index under given uncertainties, indicating it was the least resilient. In contrast, the MARL-based strategy demonstrated the lowest R-index, suggesting it was the most resilient system.

Moreover, the slopes of the contour lines varied in steepness across strategies, reflecting their relative sensitivity to noise and delay. Steeper slopes indicate that the index changes more rapidly with variations in noise, whereas gentler slopes represents greater sensitivity to delay. For example, the steep slopes observed in the centralized strategy imply that a more robust to variations in delay compared to noise. Conversely, in the decentralized strategy, if the delay is less than 15 minutes, the slopes are gentler, indicating that an increase in delay has greater negative impacts on resilience compared to noise. This suggests that the decentralized strategy is less affected by noise than the centralized strategy, as local operators independently manage their respective service areas, making disruptions easier to handle. However, the influence of noise grows as the delay extends beyond 15 minutes. The MARL-based strategy demonstrated the lowest index, indicating it is the most resilient system among the strategies. Since RH vehicles under the MARL-based strategy independently decide on relocations based on both global and local observations, it shares similarities with the decentralized strategy. As a result, the slopes of the contour lines remain gentle when the delay is less than 15 minutes.

\subsection{Efficiency of RH Services in Multi-Modal Transportation Systems During Disruptions} \label{subsection 4.5}
\begin{figure}[h]
    \centering
    \includegraphics[width=1.0\textwidth]{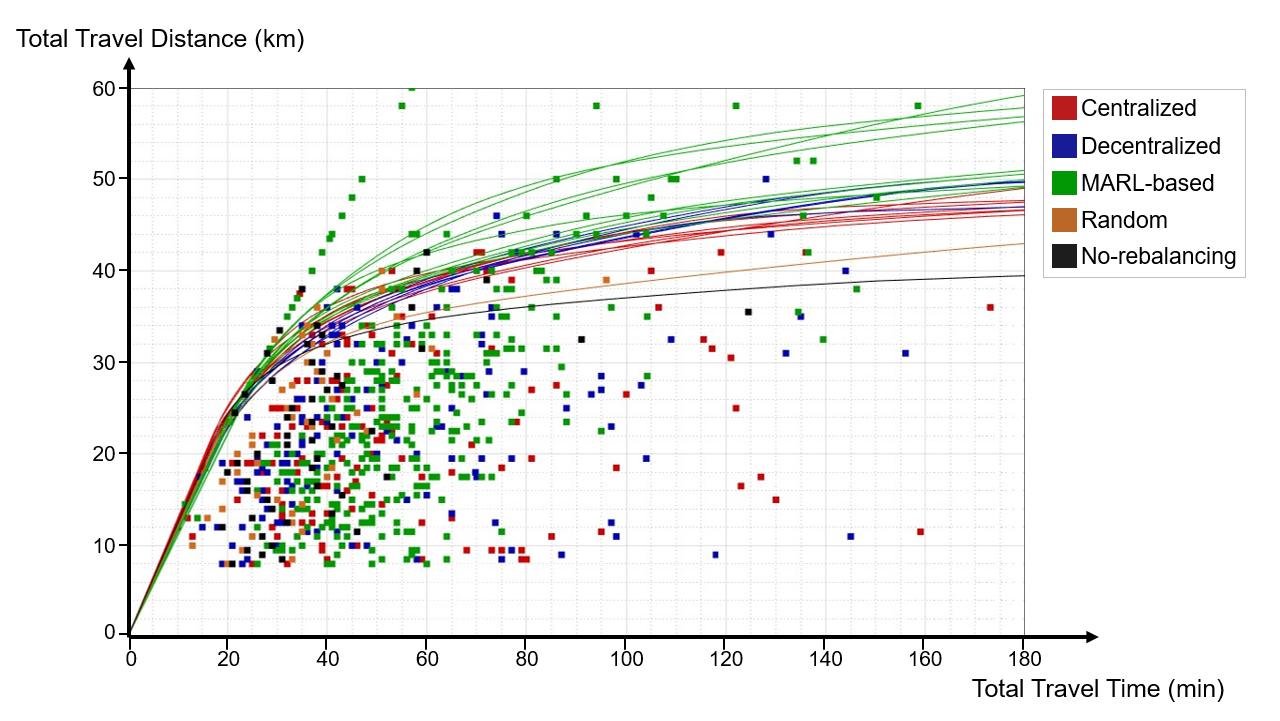}
    \caption{Comparison of total travel time and distance across rebalancing strategies in multi-modal transportation systems under disruption scenarios. The scattered dots represent user trip results with uncertainties in demand prediction noise and delay. The curves outlining the upper bounds of the dots highlight the most efficient multi-modal user trips estimated.}
    \label{Fig:Results_Travel_Time_and_Distance}
\end{figure}

RH mobility services offer a promising solution, playing a critical role in enhancing first- and last-mile connectivity. The efficiency of RH services in multi-modal transportation systems during disruption can be evaluated by analyzing user trip information. Fig.~\ref{Fig:Results_Travel_Time_and_Distance} illustrates the relationship between total travel time and total travel distance within multi-modal transportation systems under disruption scenarios.

As illustrated, the no-rebalancing strategy, represented by the black curve, showed the lowest values among the strategies. This indicates it has the least capacity for accommodating user trips. The slope of the curve reflects an efficient speed of 60 km/h for short distance trips within 25 km. However, it is important to note that these curves are based on the fastest trips, while the majority of users with the scattered dots experienced longer travel times for the same distances. The slope of the curve decreased exponentially and became almost flat as travel distances increased. This implies that it is challenging for users to travel long distances under disruption scenarios. For example, a user affected by a disruption required approximately 120 minutes for a 35 km trip. The random rebalancing strategy showed a slight improvement in operational efficiency, indicating that implementing a rebalancing strategy is more beneficial than having no strategy at all.

The MARL-based strategy, represented by green curves, showed the most efficient user trip patterns among the strategies. This strategy effectively delivered the highest number of users at an efficient speed of 60 km/h. Notably, it enabled users to travel a distance of 40 km in approximately 60 minutes, which is 3 times faster than the no-rebalancing strategy. The decentralized strategy with blue curves was assessed as the next most effective strategy. While it outperformed the MARL-based strategy in certain scenarios with prediction noise and delay, this superiority was not consistent across all conditions. The centralized strategy with red curves followed as the subsequent most effective strategy.

These findings are consistent with the results of the previous sections regarding system performance and R-indices. To enhance the resilience of multi-modal transportation systems, it is essential to implement RH rebalancing strategy effectively. More importantly, improving the accuracy of demand prediction is critical. Discrepancies between actual and predicted demand can significantly deteriorate operational efficiency, leading it difficult for even advanced rebalancing strategies to manage effectively.

\section{Conclusion} \label{section5} 

This study evaluates how RH rebalancing strategies enhance the resilience of multi-modal transportation systems during disruptions on train lines. We proposed a MARL-based RH rebalancing strategy and it was compared with various other strategies - centralized, decentralized, random, and no-rebalancing.

The contributions of the findings are as follows.
RH mobility services are incorporated into multi-modal transportation systems where users transfer between transportation modes based on their preferences. To assess RH system operations, uncertainties in demand prediction, including variations in noise and delays in response to disruptions, are incorporated to evaluate system performance during disruptions.
After demonstrating the reliability of the model in terms of traffic dynamics, system resilience is assessed through specific key performance indicators from both operator and user perspectives. By tracking vehicle and user trajectories, metrics such as user waiting time, resilience indicators, total travel time, and total travel distance are measured.
Considering that the system performance largely depends on RH rebalancing strategies, a MARL-based strategy with the MADDPG method is proposed to capture the stochastic supply-demand dynamics in large-scale networks. To minimize vehicle vacant time, a higher positive reward is granted for quicker matches. Conversely, a greater negative penalty is applied as vacant time extends. Utilizing global and local observations, the model enables agents to communicate and coordinate for efficient operations. The results demonstrate significant improvements in key performance indicators from both the operator and the user perspectives.

There are several issues that merit future investigation. Future work should focus on reducing uncertainties in demand prediction, particularly noise and delays, to further enhance the robustness of RH rebalancing systems. Additionally, improvements in computational methods or optimization algorithms could be explored to the creation of more precise and robust systems. Finally, this study was conducted using a simulated toy network, which may limit the generalizability of its findings. Future research could improve the study's validity by incorporating real-world data from diverse multi-modal transportation systems.

\vspace{0.7cm}
\appendix
\section{Appendix: Centralized System}
\label{Appendix:Centralized} 
For the centralized strategy, a main operator operates RH services for the whole network. The objective function is mathematically expressed as follows.\\

\noindent Objective function
\begin{equation}
\text{Maximize } Z = \sum_{r \in R} \sum_{v \in V} b_{rv} (U_o (r) + U_v (r)) \label{eq:Z}
\end{equation}

\noindent subject to,
\begin{flalign}
    & \sum_{v \in V} b_{rv} \leq n_r, \quad \forall r \in R \label{eq:Z_sub1} \\
    & \sum_{r \in R} b_{rv} \leq 1, \quad \forall v \in V \label{eq:Z_sub2} \\
    & b_{rv} \in \{0, 1\} \label{eq:Z_sub3}
\end{flalign}

\noindent where,  \\ 
\hspace*{0.85cm} \( U_o (r) \): operator utility for service area \( r \) \((r \in R)\)    \\ 
\hspace*{0.85cm} \( U_v (r) \): utility of RH vehicle \( v \) for service area \( r \)    \\ 
\hspace*{0.85cm} \( b_{rv} \): a binary variable indicating whether area \( r \) published a rebalance offer to vehicle \( v \)    \\ 
\hspace*{0.85cm} \( n_r \): vacancy of service area \( r \)    \\ 

The objective function in Equation~\ref{eq:Z} maximizes the total sum of operator and driver utilities. Equation~\ref{eq:Z_sub1} imposes a constraint on demand fulfillment, ensuring that the total number of vehicles that receive a rebalancing offer from service area \( r \) is equal to or less than its available capacity. Equation~\ref{eq:Z_sub2} limits assignments so that each vehicle receives at most one rebalance offer from the service areas. \( b_{rv} \) in Equation~\ref{eq:Z_sub3} is a binary variable indicating whether area \( r \) has issued a rebalance offer to vehicle \( v \). Rebalancing is performed by identifying the optimal matching of vacant vehicles to depots that yield the highest system utility \citep{duan2020centralized}.

\section{Appendix: Decentralized System}
\label{Appendix:Decentralized}
The decentralized strategy involves several local operators, each managing a separate service area. They negotiate with drivers to maximize their objectives. It is unlikely that any individual operator or driver will adhere to matching assignments based on social costs. They are primarily motivated by their self-interest and aim to maximize individual utilities \citep{chau2020decentralized}. An auction-based Gale-Shapley matching method is adopted to optimize the rebalancing problem. This two-sided stable matching algorithm has been practically applied in a variety of real-world situations (\citealp{gale1962college}; \citealp{brito2006distributed}; \citealp{ostrovsky2015fast}).

The rebalancing matching process in Algorithm~\ref{algo:Decentral} consists of four steps. The first step is setting the initial value. For operators, a set of vacancies for depots, \( N \), is obtained, with service area \( r \) having vacancy \( n_r \). The vacancy does not exceed the capacity \( n_r^o \), while considering vehicles already parked in the depots. For drivers, a set of idle vehicles, \( V_{\text{idle}} \), is obtained, consisting of vehicle \( v_i \in V \). The second step is matching operators and vehicles for rebalancing. Preferential orders of operators and drivers are calculated based on utility functions as \( \mathcal{O}_{r} \) and \( \mathcal{O}_{v_i} \), respectively. For each vehicle \( v \in V \), a preferential order \( \succ_v \) is defined for all possible matching options with service area \( r \in R \). For example, “\( (v_1, r_1) \succ_{v_1} (v_1, r_2) \succ_{v_1} (v_1, r_3) \succ_{v_1} \dots \)” indicates that vehicle \( v_1 \) prefers service area \( r_1 \) the most, then followed by \( r_2 \) and \( r_3 \), and so on. Similarly, for the service area, “\( (r_1, v_1) \succ_{r_1} (r_1, v_2) \succ_{r_1} (r_1, v_3) \succ_{r_1} \dots \)” is applied in the same process. The matching list of service area \( r \), \( M_r \), is set empty. With the preference lists given, for each \( k \)-th round (\( k = 1, \dots, K \)), individual vehicle \( i \) bids to the most preferred operator, for example \( r \). Then the vehicle \( i \) is included in the matching list \( M_r \), and the preferred element \( (r, v_i) \) is excluded from the preferential orders for vehicle \( i \). Afterwards, for operators, if the number of elements in the matching list exceeds its vacancy, candidate vehicles are excluded from the bottom rank. This process is repeated until meeting the convergence criterion or predefined iteration number. The third step is finalizing the matching between the operators and the vehicles when the convergence condition or the iteration is satisfied. Lastly, in the fourth step, the number of vehicles for rebalancing to each area is decided. The overall process is mathematically described as below.
\FloatBarrier
\renewcommand\algorithmicfor{for}
\renewcommand\algorithmicendfor{end for}
\renewcommand\algorithmicwhile{while}
\renewcommand\algorithmicendwhile{end while}
\renewcommand\algorithmicdo{}
\begin{algorithm}
\caption{Decentralized RH Rebalancing Strategy \hfill Description} \label{algo:Decentral}
\begin{algorithmic}[1]
\STATE \textbf{Step 1: Setting initial value of operators and drivers.}
\STATE $N = \{n_r \mid n_r \leq n^o_r, \; r \in R\}$ \hfill - - - - - - - - - - - - - - - - - - - - - - - - - - - - - - - - - - - - - - - obtain vacancy of depots
\STATE $V_{\text{idle}} = \{v_i \mid v_i \in V\}$ \hfill - - - - - - - - - - - - - - - - - - - - - - - - - - - - - - - - - - - - - - - - - - - - obtain idle RH vehicles

\STATE \textbf{Step 2: Matching operators and vehicles for rebalancing.}
\STATE $\mathcal{O}_r = \{(r_r, v_i) \mid (r_r, v_i) > r_r, \; (r_r, v_j)\}, \; \forall i, j \in V$ \hfill - - - - - - - - - - - - - - - - - obtain preferential orders for operators
\STATE $\mathcal{O}_{v_i} = \{(v_i, r_j) \mid (v_i, r_j) > v_i, \; (v_i, r_l), c \in V, j, l \in R\}$ \hfill - - - - - - - - - - - - - - - obtain preferential orders for drivers
\STATE $M_r = \emptyset, \; \forall r \in R$ \hfill - - - - - - - - - - - - - - - - - - - - - - - - - - - - - set a matching list of service area $r$ as an empty set

\FOR{$k = 1 : K$ \textbf{or} Meanwhile convergence not met}
    \FOR{$i = 1 : V$} 
        \STATE $r = \{r_r \mid \operatorname{rank}(r_r, O_{v_i}) = \min (\operatorname{rank}(r_r, O_{v_i}))\}, \; \forall r_r \in R$ \hfill - - - - - - - - - - find vehicle $i$'s most preferred area $r$
        \STATE $M_r \leftarrow M_r \cup v_i$ \hfill - - - - - - - - - - - - - - - - - - - - - - - - - - add vehicle $i$ into the matching list of service area $r$ 
        \STATE $\mathcal{O}_{v_i} \leftarrow \mathcal{O}_{v_i} - (r, v_i)$ \hfill - - - - - - - - - - - - - - - - - - - - - - update vehicle $i$'s preferential order by removing area $r$
    \ENDFOR
    \FOR{$r = 1 : R$}
        \WHILE{$n(M_r) > n_r$}
            \STATE $v_l = \{v_i \mid \operatorname{rank}(v_l, O_r) = \min (\operatorname{rank}(v_m, O_r))\}, \; \forall v_m \in M_r$ \hfill - - - - - - find area $r$'s most preferred vehicle $l$
            \STATE $M_r \leftarrow M_r - v_l$ \hfill - - - - - - - - - - - - - - - - - - - - - - update area $r$'s preferential order by removing vehicle $l$
        \ENDWHILE
    \ENDFOR
\ENDFOR

\STATE \textbf{Step 3: The algorithm stops when the iteration limit is reached or the convergence condition is satisfied.}

\STATE \textbf{Step 4: Determine the number of rebalancing vehicles by service areas.}

\end{algorithmic}
\end{algorithm}
\FloatBarrier

\vspace{0.7cm}
\noindent {\large \textbf{CRediT authorship contribution statement}} \\
\noindent \textbf{Euntak Lee}: Conceptualization, Investigation, Methodology, Formal analysis, Visualization, Validation, Writing - original draft, Writing - review \& editing. \textbf{Rim Slama}: Conceptualization, Methodology, Writing - review \& editing. \textbf{Ludovic Leclercq}: Conceptualization, Methodology, Formal analysis, Investigation, Validation, Supervision, Writing - review \& editing, Project administration, Funding acquisition.

\vspace{0.7cm}
\noindent {\large \textbf{Declarations of interest}} \\
\noindent The authors declare that they have no known competing financial interests or personal relationships that could have appeared to influence the work reported in this paper.

\vspace{0.7cm}
\noindent {\large \textbf{Acknowledgment}} \\
\noindent This research has received funding from the European Union’s Horizon Europe research and innovation program under Grant Agreement no. 101103808 (ACUMEN).

\bibliographystyle{cas-model2-names}
\bibliography{references} 

\end{document}